\documentclass[english,11pt]{article}

\usepackage[german,french,english]{babel}
\usepackage[cp850]{inputenc}
\usepackage{latexsym,graphicx, fancybox}
\usepackage{graphicx}
\usepackage{amssymb, amsmath, amsfonts}
\usepackage{color}
\usepackage{latexsym}

\font\tenmath=msbm10 scaled 1200
\font\sevenmath=msbm7 scaled 1200
\font\fivemath=msbm5 scaled 1200 

\newfam\mathfam \textfont\mathfam=\tenmath
\scriptfont\mathfam=\sevenmath \scriptscriptfont\mathfam=\fivemath

\def\R{{\mathbb R}}

\def\E{{\mathbb E}}

\def\P{{\mathbb P}}

\def\C{{\mathbb C}}

\def\F{{\cal F}}

\newtheorem{Theorem}{Theorem}[section]
\newtheorem{Definition}[Theorem]{Definition}
\newtheorem{Proposition}[Theorem]{Proposition}

\newtheorem{Lemma}[Theorem]{Lemma}

\newfam\mathfam \textfont\mathfam=\tenmath
\scriptfont\mathfam=\sevenmath \scriptscriptfont\mathfam=\fivemath

\def \^#1{\if#1i{\accent"5E\i}\else{\accent"5E#1}\fi}

\def \a{\alpha}

\def \s{\sigma}

\begin{document}
\selectlanguage{english}
\title{\bf  Volterra equations with   affine drift: looking for stationarity. Application to   the quadratic rough Heston model}
 
\author{
{\sc  Gilles Pag\`es} \thanks{Sorbonne Universit\'e, Laboratoire de Probabilit\'es, Statistique et Mod\'elisation, case 158, 4, pl. Jussieu, F-75252 Paris Cedex 5, France. E-mail: {\tt  gilles.pages@sorbonne-universite.fr}}~\thanks{This research
benefited from the support of the ``Chaire Risques Financiers'', Fondation du Risque.}  
}
\maketitle 
\renewcommand{\abstractname}{Abstract}
\begin{abstract}
We investigate the properties of  the solutions   of scaled Volterra equations (i.e. with an affine mean-reverting drift) in terms of  stationarity at both a finite horizon and on the long run. In particular we prove that such an equation never has a  stationary regime, except if the  kernel is constant  (i.e. the equation is a standard Brownian diffusion) or  in some fully degenerate pathological settings. We introduce a deterministic stabilizer $  \varsigma$ associated to the kernel which may produce a {\em fake stationary regime} in the sense that all the marginals share  the same expectation and variance. We also show that  the marginals  of such a process when starting various initial values  are confluent in $L^2$ as time goes to infinity. We   establish that for some classes of diffusion coefficients (square root of positive quadratic polynomials)  the time shifted solutions of such Volterra equations weakly functionally converge toward a family of  $L^2$-stationary  processes sharing the same covariance function. We apply these results to  a  family  of stabilized quadratic rough  volatility models (when the kernel $K(t)= \frac{t^{H-\frac 12}}{\Gamma(H+\frac 12)}$, $0<H<\frac 12$) which  produces   examples Volterra processes sharing a fake stationary regime.
\end{abstract}
\section*{Introduction}  
We  investigate in this paper  conditions under which   stochastic Volterra equations with affine drifts and convolutive kernels of the form
\[
X_t= X_0 +\int_0^t K(t-s)(\mu(s)-\lambda X_s)ds+ \int_0^t K(t-s)\sigma(X_s)dW_s
\]
have a stationarity regime  in the classical sense, i.e. the distribution of the process is invariant by time shift (see~\cite{Jacquieretal2022}) or in a weaker sense~--~that we called   {\em fake stationary regime (of type~I)}~--~in which  the solution only has constant mean and variance at every instant $t$. In the Gaussian case (typically pseudo-Ornstein-Uhlenbeck process),  it implies that the process has the same $1$-marginal distribution at very time $t$ (which could be called a {\em fake stationary regime of type~II}). However such a property if it holds does not imply stationarity nor classical weak $L^2$-stationarity (based on covariance).

In terms of application this may have consequences on rough volatility models (see~\cite{BayerFG2016,ElEuchR2018,GatheralJR2018,GaJuRo2020}) which met a great success in the financial community during these last ten years, in particular because it provides a consistent modeling of the dynamics of financial  markets (mostly equity) at various scales from the order book to the pricing of derivative products in spite of its non-Markovian feature. In the  original version of these  {\em rough} stochastic volatility models, mostly derived from the regular Heston model, the  mean-reverting CIR dynamics of the volatility is replaced by a Volterra version of this CIR process in which the convolution kernel is singular and of  fractional  form $K(t)=t^{H-\frac 12}$ for some $H\!\in (0,1/2)$ (Hurst coefficients). In the late 1990's, attempts to model the long memory effect on continuous time stochastic volatility models (see~\cite{ComteR1996, ComteR1998}) also relied on similar stochastic Volterra     equations using the same family of kernels but with $H\!\in (1/2,1)$. In both cases the aim is to mimic the behaviour of a stochastic differential equations (SDE) driven by a fractional Brownian motion without having to face the technicalities of such  equations which rely on rough path theory, technicalities which can become extreme when $H$ gets close to $0$.

If we come  back for  a while to the standard Heston model equipped with its CIR driven stochastic volatility process,  it has been observed  for long by practitioners  that this model was difficult to calibrate (on vanilla products)  to produce good prices for derivatives with  short  maturities, especially path-dependent ones. This is most likely due to the fact that the CIR process is  usually initialized by a deterministic value at the origin (often the  long run mean). So this volatility dynamics  actually has two regimes, first, one for short maturities when the starting value is prominent and the variance  is very small, then it deploys until attaining its  its stationary regime~--~ a $\gamma$-distribution~--. A solution  to get rid of this drawback is to consider directly a Heston model in which the volatility is taken  under its stationary regime (see e.g.~\cite{PagPan2009}) which also has the advantage to reduce by one the number of parameters. A quantization based approach is provided for calibration and pricing of path-dependent derivatives  (at least under the Feller condition that ensures positivity of the volatility) in a recent paper~\cite{LemMonPag2022}.

If, as it is likely,  a similar problem occurs for rough volatility models based on the  rough  Volterra version of the CIR model,  the search for some  stationarity properties of such a rough CIR model can be of interest to make possible a calibration  consistent  for both short and long maturities. This does not provide a plug-and-play solution to this problem but explores in an elementary way  what are  or could be the connexions between stationarity and Volterra equations with  ``mean-reverting'' affine drift with a mitigate conclusion. This could be compared to the recent paper~\cite{Jacquieretal2022} about ergodicity of Volterra processes. 

\medskip
The paper is organized as follows. In Section~\ref{sec:background} we recall the main properties of  stochastic Volterra equations with convolutive drifts: existence, moment control, etc, with a focus on those of these processes having an affine drift. In that case some  specific tools are available like the Resolvent and the solution  Wiener-Hopf equation which corresponds when $K(t)= e^{\lambda t}$ (Markovian framework) to apply It\^o's lemma to $e^{\lambda t}X_t$. In Section~\ref{sec:genfakestatio} we elucidate when a Volterra equation with affine  has an invariant distribution  (distribution invariant by time shift). We show  that this situation is essentially degenerate in the sense that it  implies that the kernel $K$ is constant i.e. that the Volterra equation is a standard Brownian equation. In Section~~\ref{sec:TowardFake} we define some notions of  {\em fake} stationarity   like having constant mean and variance (type~I), or constant $1$-marginal distribution (type~II). To this end we are led to introduce a multiplicative deterministic function$\varsigma$  in the Brownian convolution, called {\em stabilizer} to ``control" the process (through its volatility coefficient). This stabilizer is solution to an intrinsic convolution equation related to the kernel of the equation. In  Section~\ref{sec:FakeI-II} we provide examples  of  such``fake stationary'' stabilized processes with some explicit computations  when 
\[
\sigma(x) = \sqrt{ \kappa_0 +\kappa_1\,(x-\tfrac{\mu_0}{\lambda})+\kappa_2\,(x-\tfrac{\mu_0}{\lambda})^2}\quad \mbox{ with }\quad  \kappa_i\ge 0,\;i=1,2, \;\kappa^2_1 \le 4\kappa_0\kappa_2. 
\]
This form of the diffusion coefficient appears e.g. in~\cite{GaJuRo2020} when modeling the volatility in the quadratic rough Heston model. We also  investigate for such stabilized processes the  functional weak asymptotics of the time shifted process $(X_{t+s})_{s\ge 0}$ as $t\to +\infty$ which turns out to be a weak $L^2$-stationary process. Finally in Section~\ref{sec:Back2ML} we apply our results to the case case of   fractional kernels $K_{\alpha}(t)= t^{\alpha-1}$, $\alpha\, \in (0,1)$ covering  ``rough''  models (where $\a = \frac 12 +H$, $H\!\in (0,1/2)$, the Hurst coefficient). In particular we prove the existence of the stabilizer with a semi-closed form as (fractional) power series whose coefficients can be computed by induction. We prove in an appendix that this series has an infinite convergence radius, null at $0$, increasing toward a finite (identified) limit at infinity.

 \medskip
 \noindent {\sc Notations.} 
$\bullet$ ${\rm Leb}_1$ (resp. ${\rm Leb}_{\R_+}$) denotes the Lebesgue measure on $\R$ (resp. $\R_+= [0, +\infty)$). 

\smallskip 
\noindent $\bullet$ For $f, g \!\in {\cal L}_{\R_+,loc}^1 (\R_+, {\rm Leb}_1)$, we define their convolution by $f*g(t) = \int_0^tf(t-s)g(s)ds$, $t\ge 0$.

\smallskip 
\noindent $\bullet$ For $f, g \!\in {\cal L}_{\R_+,loc}^2(\R_+, {\rm Leb}_1)$ and $W$ a Brownian motion, we define  their stochastic convolution by 
\[
f\stackrel{W}{*}g = \int_0^t f(t-s)g(s) dW_s, \quad t\ge 0.
\]
 
\noindent $\bullet$ $[X]$ stands for the distribution of the  random variable/vector/process $X$.

\smallskip 
\noindent $\bullet$ $X\perp \! \! \!\perp Y$   stands for the independence of the two random objects $X$ and $Y$.  

\section{Background on   Volterra SDEs with convolutive kernels}\label{sec:background}
\subsection{Volterra processes with convolutive kernels}

We are interested in the stochastic Volterra equation
 \begin{equation}\label{eq:Volterra}
   X_t=X_0+\int_0^t K(t-s) b(s,X_s)ds+\int_0^t K(t-s)\sigma(s,X_s)dW_s,\quad t\ge 0,
 \end{equation}
where $b:[0,T]\times\R\to \R$, $\sigma:[0,T]\times\R\to\R$ are Borel  measurable, $K\!\in {\cal L}^2_{loc,\R_+}({\rm Leb}_1)$  is a convolutive kernel and $(W_t)_{t\ge 0}$ is a standard  Brownian motion independent from the $\R$-valued random variable $X_0$ both defined on a probability space $(\Omega,{\cal A}, \P)$. Let $(\F_t)_{t\ge 0}$ be a filtration (satisfying the usual conditions) such that $X_0$ is ${\cal F}_0$-measurable  and $W$ is an $({\cal F}_t)$-Brownian motion independent of $X_0$.

Assume that  the kernel $K$ satisfies   
\begin{equation}\label{eq:Kint}
\exists\, \beta >1\quad \mbox{such that}  \quad K\!\in L^{2\beta}_{loc}({\rm Leb}_1)
\end{equation}
and that,   for every $T>0$,
\begin{equation}\label{eq:Kcont}
({\cal K}^{cont}_{T,\theta}) \;\exists\, \kappa_{_T}< +\infty,\; \exists\, \theta_{_T} >0,\; \forall\,\delta \!\in (0,T),\;
 \sup_{t\in [0,T]} \left[\int_0^t |K(\big(s+\delta)\wedge T\big)-K(s))|^2ds \right]^{\frac 12} \le  \kappa_{_T}\,\delta^{\theta_{_T}}.
\end{equation}

Assume $b$ and $\sigma$ satisfy the following Lipschitz-linear growth assumption uniform in time
\begin{align*}
(i) &\; \forall\, t\!\in [0,T], \;  \forall\, x,\, y \!\in \R, \; |b(t,x)-b(t,y)+|\sigma(t,x)-\sigma(t,y)| \le C_{b,\sigma,T}|x-y|,\\
(ii)& \sup_{t\in [0,T]}|b(t,0)| +|\sigma(t,0)| <+\infty.
\end{align*}

It has been established in~\cite[Theorem~1.1]{JouPag22}  that, if $X_0\!\in L^p(\P)$ for some $p>0$, then  Equation~\eqref{eq:Volterra}  admits  a unique pathwise continuous solution on $\R_+$  starting from   $X_0$ satisfying (among other properties),
\begin{equation}\label{eq:L^p-supBound}
\forall\, T>0, \; \exists \,C_{_{p,T} }>0,\quad \Big\| \sup_{t\in[0,T]}|X_t|\Big\|_p \le C_{_{p,T} }(1+\|X_0\|_p).
\end{equation}
This result appears as a generalization of the classical  strong existence-uniqueness result of pathwise continuous solutions   established in~\cite{ZhangXi2010}, only when the starting value $X_0$ has finite polynomial moments at any order (The framework is more general).
\subsection{Laplace Transform}
Let us first introduce the Laplace transform of a Borel function $f:\R_+\to \R_+$ by
\[
\forall\, t\ge  0, \quad L_f(t)= \int_0^{+\infty} e^{-tu}f(u)du \!\in [0,\infty].
\] 
This Laplace transform is non-increasing and if $L_f(t_0)<+\infty$ for some $t_0\ge 0$, then $L_f(t)\to 0$ as $t\to +\infty$.
 
 One can define the  Laplace transform of a Borel function $f:\R_+\to \R$ on $(0, +\infty)$ as soon as $L_{|f|}(t) <+\infty$ for every $t>0$  by the above formula. The Laplace transform can be extended to $\R_+$ as an $\R$-valued function  if $f \!\in {\cal L}^1_{\R_+}({\rm Leb}_1)$. 

\subsection{Resolvent of a convolutive kernel}
For every $\lambda \!\in \R$,  the {\em resolvent} $R_{\lambda}$ associated to $K$ and $\lambda$ is defined as the unique solution -- if it exists --  to the deterministic Volterra equation
\begin{equation}\label{eq:Resolvent}
\forall\,  t\ge 0,\quad R_{\lambda}(t) + \lambda \int_0^t K(t-s)R_{\lambda}(s)ds = 1
\end{equation}
or, equivalently, written in terms of convolution, 
\begin{equation}\label{eq:Resolvent2}
R_{\lambda}+\lambda K*R_{\lambda} = 1.
\end{equation}
The solution  always satisfies $R_{\lambda}(0)=1$ and formally reads 
\begin{equation}\label{eq:Resolvent3}
R_{\lambda} = \sum_{k\ge 0} (-1)^k \lambda^k (\mbox{\bf 1}*K^{*k})
\end{equation}
where, by convention in this formula, $K^{*0}= \delta_0$ (Dirac mass at $0$).

\medskip From now on we will assume that the kernel $K$ satisfies
\begin{equation}\label{eq:LaplaceK}
({\cal L}_K)\qquad\forall \, t>0, \quad L_K(t)<+\infty.
\end{equation}

\bigskip
\noindent  $\bullet$   If $\lambda >0$,~\eqref{eq:Resolvent3} defines an absolutely convergent  series for every $t\!\in (t_K, +\infty)$  where $t_K=\inf\{t: L_K(t)\le1/\lambda\}<+\infty$.

\medskip
\noindent $\bullet$  If the (non-negative) kernel $K$  satisfies 
\begin{equation}\label{eq:Kcontrol}
0\le K(t)\le Ce^{bt }t^{a-1}\mbox{  for some }\; a,\, \; b,\; C>0\!\in \R_+ 
\end{equation}then one easily checks by induction  that  
\begin{equation}\label{eq:BoundK}
 \mbox{\bf 1}*K^{*n}(t) \le C^n e^{bt}\frac{\Gamma(a)^n}{\Gamma(an+1)} t^{an},
\end{equation}
 where $\Gamma(a)= \int_0^{+\infty}e^{-u}u^{a-1}du$. Then, one shows  using Stirling formula    that, for such kernels, the above series~\eqref{eq:Resolvent3} is absolutely converging for every $t>0$  so that the function   $R_{\lambda}$ is well-defined on $(0, +\infty)$. 
 
 \noindent $\bullet$ If $K$ is regular enough (say continuous) the resolvent $R_{\lambda}$ is differentiable and one checks that $f_{\lambda}=-R'_{\lambda}$ satisfies, for every  $t>0$, 
 \[
-f_{\lambda}(t) +\lambda \big( R_{\lambda}(0)K(t) - K *f_{\lambda}(t)\big)=0
\]
that is $f_{\lambda}$ is solution to the equation
\begin{equation}\label{eq:flambda-eq}
f_{\lambda} +\lambda K *f_{\lambda}=\lambda   K.
\end{equation}
In particular, if $R_{\lambda}$ turns out to be non-increasing, then $f_{\lambda}$ is non-negative and satisfies $0\le f_{\lambda} \le \lambda K$. In that case one also has that $\int_0^{+\infty} f_{\lambda}(t)dt = 1 -R_{\lambda}(+\infty)$, so that $f_{\lambda} $ {\em is a probability density} if and only if $\displaystyle \lim_{t\to +\infty} R_{\lambda}(t) =0$.

\bigskip
 \noindent {\bf Examples of kernels}. 
\begin{enumerate}

\item {\em Trivial kernel (Markov)}. Let $K(t) = \mbox{\bf 1}_{\R_+}(t)$. It obviously satisfies~\eqref{eq:Kint},~\eqref{eq:Kcont} and~\eqref{eq:Kcontrol}. Then $R_{\lambda} (t)= e^{-\lambda t}$. 
\item  {\em Exponential kernel}.  Let $K(t)= e^{-u} $.  It clearly satisfies~\eqref{eq:Kint}, ~\eqref{eq:Kcont} and~\eqref{eq:Kcontrol}. Then $R_{\lambda} (t)= \left\{\begin{array}{cc}t+1& \mbox{if }\lambda =-1\\ \frac{1+\lambda e^{-(\lambda+1)t}}{\lambda +1} & \mbox{if }\lambda \neq -1 \end{array}\right.$.
\item {\em Fractional integration kernel}. Let 
\begin{equation}\label{eq:frackernel}
K(t) = K_{\alpha}(t) = \frac{u^{\alpha-1}}{\Gamma(\alpha)} \mbox{\bf 1}_{\R_+}(t),  \quad \alpha>0.
\end{equation}
   These kernels  satisfy~\eqref{eq:Kint} and~\eqref{eq:Kcont} for $\alpha >1/2$ (with $\theta_{_T}= (\alpha-\frac 12)\wedge 1$, see~\cite{RiTaYa2020} or~\cite{JouPag22} along many others)~among many others and trivially~\eqref{eq:Kcontrol}. This  family of  kernels corresponds to the fractional integrations of order $\alpha >0$. 

It follows from the  easy identity  $K_{\alpha}*K_{\alpha'}= K_{\alpha+\alpha'}$ that 
\begin{equation}\label{eq:e_alpha}
R_{\alpha,\lambda}(t) = \sum_{k\ge 0} (-1)^k \frac{\lambda ^k t^{\alpha k}}{\Gamma(\alpha k+1)}= E_{\alpha}(-\lambda t^{\alpha} ) \; t\ge 0,
\end{equation}
where $E_{\alpha}$ denotes the {\em standard   Mittag-Leffler function}  
\[
E_{\alpha}(t) = \sum_{k\ge 0} \frac{t^k}{\Gamma(\alpha k+1)},\ t\!\in \R.
\] 
One shows (see~\cite{GorMain1997} and Section~\ref{sec:Back2ML} further on) that $E_{\alpha}$ is increasing and differentiable on the real line with $\displaystyle \lim_{t\to +\infty}E_{\alpha} (t) =+\infty$ and $E_{\alpha}(0)=1$. In particular, $E_\alpha$ is an homeomorphism from $(-\infty, 0]$ to $(0,1]$. Hence, if $\lambda >0$, the function $f_{\alpha, \lambda}$ defined on $(0,+\infty)$ by
\begin{align}
\label{eq:flambda} f_{\alpha, \lambda}(t) &= - R'_{\alpha, \lambda}(t) = \alpha\lambda t^{\alpha-1} E'_{\alpha}(-\lambda t^{\alpha})  = \lambda t^{\alpha-1}\sum_{k\ge 0}(-1)^k\lambda^k \frac{t^{\alpha k}}{\Gamma(\alpha (k+1))}
\end{align}
is a probability density -- called Mittag-Leffler density -- since $f_{\alpha, \lambda}>0$ and $\displaystyle \int_0^{+\infty} f_{\alpha, \lambda}(t)dt =  R_{\alpha,\lambda}(0) - R_{\alpha,\lambda}(+\infty) = 1$. Note that when $\alpha =1$ (i.e. $K=K_1=\mbox{\bf 1}$), $E_{1}(t) = e^t$, $R_{1,\lambda}(t)= e^{-\lambda t}$ and $f_{1,\lambda}(t)=\lambda e^{-\lambda t}$.
%
\item {\em Exponential-Fractional integration kernel}. Let $K(t) = K_{\alpha,\rho}(t) = e^{-\rho t} \frac{u^{\alpha-1}}{\Gamma(\alpha)} \mbox{\bf 1}_{\R_+}(t)$ with $\alpha, \rho>0 >0$. One checks that these kernels also satisfy~\eqref{eq:Kint} and~\eqref{eq:Kcont} for $\alpha >1/2$ (still with $\theta_{_T}= (\alpha-\frac 12)\wedge 1$) and trivially~\eqref{eq:Kcontrol}. Moreover 
$K_{\alpha, \rho}*K_{\alpha', \rho}= K_{\alpha+\alpha',\rho}$ so that the resolvent reads
\[
R_{\alpha, \rho, \lambda}(t)= e^{-\rho  t}R_{\alpha, 0,\lambda}(t).
\]
\end{enumerate}

\medskip
This last example is important for applications, especially in Finance for  the wide class of rough   stochastic volatility models (see~\cite{JaissonR2016,GatheralJR2018,ElEuchR2018,GaJuRo2020}) since in that case $\alpha =H+\frac 12$ where $H\!\in (0,1)$ is kind of Hurst coefficient representative or the roughness of the model ($H<\frac 12$) or its long memory ($\frac 12<H<1$). 

%
%
%

\subsection{Application to  the Wiener-Hopf equation}
We come now to the main result of these preliminaries. 
\begin{Proposition}[Wiener-Hopf equation] \label{prop:W-H}  Let $g: \R_+ \to \R$ be a  locally bounded Borel function, let $K \! \in L^1_{loc}({\rm Leb}_{\R_+})$ and let $\lambda \!\in \R$.  Assume that $R_{\lambda}$ is differentiable on $(0, +\infty)$ with a derivative $R'_{\lambda}\!\in     L^1_{loc}({\rm Leb}_1{\R_+})$, that both $R_{\lambda}$ and $R'_{\lambda}$ admit a finite Laplace transform on $\R_+$ and $\displaystyle \lim_{u\to +\infty} e^{-tu}R_{\lambda}(u) = 0$ for every $t>0$. Then, the Wiener-Hopf equation
\begin{equation}\label{eq:Wiener-Hopf}
\forall\, t\ge 0, \quad x(t) = g(t) -\lambda \int_0^t K(t-s) x(s) ds
\end{equation}
(also reading $x= g-\lambda K*x$)  has   a solution given by
\begin{equation}\label{eq:Wiener-Hopf-solu}
\forall\, t\ge 0, \quad x(t) = g(t) +\int_0^t R'_{\lambda}(t-s)g(s)ds
\end{equation}
or, equivalently, 
$$
x= g- f_{\lambda}*g,
$$
where $f_{\lambda}= -R'_{\lambda}$. This solution is uniquely defined on $\R_+$ up to  $dt$-$a.e.$ equality.
\end{Proposition}

We provide a proof of this classical result for the reader's convenience (see~e.g.~\cite{GaJuRo2020}).

\bigskip
\noindent {\bf Proof.}   First  note that, by an integration by parts,  that for every $t>0$, 
\[
tL_{R_{\lambda}} (t)= 1 +L_{R'_{\lambda}}(t). 
\]
On the other hand it follows from~\eqref{eq:Resolvent}  that $tL_{R_\lambda}(t)(1+\lambda L_K(t))=1 $, $t>0$. Consequently,
\[
L_x (t)= \frac{L_g(t)}{1+\lambda L_{K}(t)}= tL_g (t)L_{R_{\lambda}}(t)= L_g(t)  \big(1 +L_{R'_{\lambda}}(t)\big)= L_{g+ R'_{\lambda}*g}(t)
\]
which completes the proof since Laplace transform is injective. \hfill $\Box$

\section{Genuine and fake stationarity of a scaled Volterra equation}\label{sec:genfakestatio}

Throughout what follows we assume that the $\lambda$-resolvent $R_{\lambda}$ of the kernel $K$   satisfies for every $\lambda >0$
\begin{equation}\label{eq:hypoRlambda}
({\cal K})\; \left\{
\begin{array}{ll}
(i) & R_{\lambda} \;\mbox {decreases  to $0$
as $\, t\to+\infty$ and } f_{\lambda}=- R'_{\lambda} >0\; \mbox{ on }\;(0, +\infty).\\
(ii) &  f_{\lambda} \!\in {\cal L}_{loc}^2({\rm Leb}_1)
\end{array} \right.
\end{equation}
In particular, under this assumption,  the function $f_{\lambda}$ is a probability density which has subsequently a finite Laplace transform defined (at least) on $\R_+$.

\smallskip
This assumption is satisfied by our examples~1-3-4 of kernels. 

\medskip
From now we focus on  the special case of a {\em scaled} Volterra equation associated to a convolutive kernel $K:\R_+\to \R_+$ satisfying~$({\cal K})$,~\eqref{eq:Kint},~\eqref{eq:Kcont} and~\eqref{eq:Kcontrol}:
\begin{equation}\label{eq:Volterrameanrevert}
X_t = X_0 +\int_0^t K(t-s)(\mu(s)-\lambda X_s)ds + \int_0^t K(t-s)\sigma(s,X_s)dW_s, \quad X_0\perp\!\!\!\perp W,
\end{equation}
where $\lambda>0$, $\mu :\R_+\to \R$ is a  bounded Borel function (hence having   a well-defined finite Laplace transform on $(0,+\infty)$) and $\sigma: \R_+\times \R \to  \R$ is Lipschitz continuous in $x$, locally uniformly in $t\!\in [0,T]$, for every $T>0$. Note that the drift $b(t,x) = \mu(t)-\lambda  x$ is clearly Lipschitz continuous in $x$,  uniformly in $t\!\in \R$. 

Then, Equation~\eqref{eq:Volterrameanrevert} has a unique $({\cal F}^{X_0,W}_t)$-adapted solution $(X_t)_{t\ge 0}$ starting from $X_0\!\in L^0(\P)$ (Apply~\cite[Theorem~1.1]{JouPag22} on every time interval $[0,T]$, $T>0$. Then pasting   theses solutions taking advantage of the uniform in time linear growth of the drift and of $\sigma$.  

Have in mind  that under our assumptions, if $p>0$ and $\E\, |X_0|^p <+\infty$, then $\E\sup_{t\in [0,T]} |X_t|^p <C_{_T}(1+\|X_0\|_p)+\infty$ for every $T>0$  (see~\cite[Theorem~1.1]{JouPag22}).  Combined with the fact that $|\sigma(t,x)|\le C'_{_T}(1+|x|)$, $t\!\in [0,T]$,  this implies that 
$\sup_{t\in [0,T]}\|\sigma(t,X_t)\|_p <C(1+\|X_0\|_p)<+\infty$ for every $T>0$. This allows us to use in what follows without restriction both regular and stochastic Fubini's theorems.


\begin{Proposition}[Wiener-Hopf  transform] \label{prop:Wienr-HopfX} Let $\lambda >0$ and let $\mu :\R\to \R$ be a   bounded Borel function. Assume the kernel $K$ satisfies the above  assumptions~$({\cal K})$,~\eqref{eq:Kint},~\eqref{eq:Kcont} and~\eqref{eq:Kcontrol}  
Then, the solution $(X_t)_{t\ge 0}$ of the scaled Volterra equation~\eqref{eq:Volterrameanrevert} also satisfies
\begin{equation}\label{eq:Volterrameanrevert2}
X_t= X_0R_{\lambda}(t) +\frac{1}{\lambda}\int_0^t f_{\lambda}(t-s)\mu(s)ds + \frac{1}{\lambda}\int_0^t f_{\lambda}(t-s)\sigma(s,X_s)dW_s.
\end{equation}
\end{Proposition}

 \smallskip
 \noindent {\bf Proof.} We may read the above equation~\eqref{eq:Volterrameanrevert}  ``pathwise'' (i.e. $\omega$ by $\omega$) as a Wiener-Hopf equation with $x(t)= X_t(\omega)$ and 
\[
g(t) = X_0(\omega) + (\mu*K)_t + \big(K\stackrel{W}{*} \sigma(.,X_.(\omega))\big)_t.
\]

Then
\begin{align*}
X_t & =  X_0+  (\mu*K)_t + \big(K\stackrel{W}{*} \sigma(.,X_.)\big)_t + \int_0^t R'_{\lambda}(t-s)\Big[   X_0+ (\mu*K)_s+ \big(K\stackrel{W}{*} \sigma(.,X_.)\big)_s\Big] ds \\
&= X_0+ (\mu*K)_t + \big(K\stackrel{W}{*} \sigma(.,X_.)\big)_t +X_0 \int_0^t R'_{\lambda}(t-s)ds \\
& \qquad +\underbrace{ \int_0^t R'_{\lambda}(t-s)(\mu*K)_s ds }_{=: (a)} + \underbrace{ \int_0^t R'_{\lambda}(t-s) \big(K\stackrel{W}{*} \sigma(.,X_.)\big)_sds}_{=: (b)}.
\end{align*}

By commutativity and associativity (which relies on regular Fubini's theorem) of regular convolution, one has
\begin{equation}\label{eq:(a)init}
(a)  = -f_{\lambda}*(\mu*K)_t = - \big((f_{\lambda}*K)*\mu\big)_t.
\end{equation}
where $f_{\lambda}= -R'_{\lambda}$. On the other hand by differentiating Equation~\eqref{eq:Resolvent}  we know that  that  $f_{\lambda}$ satisfies Equation~\eqref{eq:flambda-eq}  that can be rewritten as 
\[
-f_{\lambda}*K = \frac{1}{\lambda}f_{\lambda} - K.
\]
Consequently, inserting this identity in~\eqref{eq:(a)init} yields
\begin{align}\label{eq:(a)}
(a) &  =  \frac{1}{\lambda} (f_{\lambda}*\mu)_t -(K*\mu)_t.
\end{align}
Using stochastic Fubini's theorem and, again, Equation~\eqref{eq:flambda-eq} recalled above  satisfied by $f_{\lambda}$, we derive that
\begin{align}\label{eq:(b)}
(b) &=   \frac{1}{\lambda} \big(f_{\lambda}\stackrel{W}{*} \sigma(\cdot, X_{\cdot})\big)_t - \big(K\stackrel{W}{*} \sigma(\cdot, X_{\cdot})\big)_t.
\end{align}
inserting~\eqref{eq:(a)} and~\eqref{eq:(b)} into~\eqref{eq:Volterrameanrevert} finally yields the announced 
\begin{align*}
  \hskip 4cm X_t &= X_0R_{\lambda}(t)+ \frac{1}{\lambda}(f_{\lambda}*\mu)_t + \frac{1}{\lambda}\big(f_{\lambda}\stackrel{W}{*} \sigma(\cdot,X_{\cdot})\big)_t . \hskip 4cm  \Box 
\end{align*}

\medskip
\noindent {\bf Remark.} When $K(t) = \mbox{\bf 1}$, the above computations correspond  to the application of It\^o's Lemma to $e^{\lambda t}X_t$. 

\subsection{Stationarity of the mean}\label{sec:moyennecste}  Before tackling the problem of the    stationary regime of the  ``scaled'' Volterra equation~\eqref{eq:Volterrameanrevert}, a first step  can be to determine when this equation has a constant mean if $X_0\!\in L^2(\P)$~(\footnote{$X_0\in L^p(\P)$ for some $p>1$ would be  sufficient by standard argument based on $BDG$ Inequality in this section.}) i.e.
$\E\, X_t = \E\, X_0$ for every $t\ge 0$. As, for every $t\!\in [0,T]$, $T>0$,  
\[
\E\, \bigg(\int_0^t f_{\lambda}(t-s)\sigma(s,X_s) dW_s\bigg)^2 = \int_0^t f^2_{\lambda}(t-s)\E\,\sigma^2(s,X_s) ds\le C_{_T}(1+\|X_0\|_2^2) \int_0^tf^2_{\lambda}(u)du<+\infty,
\]
we have, for every $t\ge 0$, 
$$
\E\,\int_0^t f_{\lambda}(t-s)\sigma(s,X_s) dW_s=0
$$ 
so that 
\[
\E\, X_t = R_{\lambda}(t) \E\, X_0  + \frac{1}{\lambda}(f_{\lambda}*\mu)_t.
\]
Using that $R_{\lambda}(t) = 1 -\int_0^t f_{\lambda}(t-s)ds $, we derive that, for every $t\ge 0$,
\[
\E\, X_t = \E\, X_0 +\int_0^t  f_{\lambda}(t-s)\Big(\frac{\mu(s)}{\lambda}-\E\, X_0  \Big)ds.
\]
Hence $\E\, X_t$ is constant if and only if
\[
\forall\, t\ge 0, \quad \int_0^t  f_{\lambda}(t-s)\Big(\frac{\mu(s)}{\lambda}-\E\, X_0  \Big)ds =0.
\]
This also reads in terms of Laplace transform
\[
L_{f_{\lambda}}\cdot L_{\frac{\mu(\cdot)}{\lambda}-\E\,X_0}= 0.
\]
As $L_{f_{\lambda}} (t)> 0$ on $(0,+\infty)$ owing to Assumption~(${\cal K})$$(ii)$ (see~\eqref{eq:hypoRlambda}), this implies that $\frac{\mu(t)}{\lambda}-\E\,X_0=0$ $dt$-$a.e.$ i.e.
\[
\mu(t) = \mu_0 \;\; dt\mbox{-}a.e. \quad \mbox{ and } \quad\E\, X_0 = \frac{\mu_0}{\lambda}.
\] 

\subsection{True scaled Volterra equation with time homogenous coefficient $\s$ are (usually) not stationary}  
%
If we assume that  $X_0\!\in L^2(\P)$, under our assumptions, $X$ is well-defined on the whole positive real  line and $X_t\!\in L^2(\P)$ for every $t\ge 0$.  In fact  (see~\eqref{eq:L^p-supBound}), $\|\sup_{t\in [0,T]}|X_t|\|_2 < C_T(1+\|X_0\|_2)<+\infty$ for every $T>0$ where $C_T$ is a (non-exploding) positive real constant. 
Then   the non-negative function defined by 
   \begin{equation}\label{eq:barsigma}
t\longmapsto \bar \sigma^2(t):= \E\,\sigma^2(t,X_t), \quad t\ge 0.
\end{equation}
is  locally bounded   on $\R_+$  since $\sigma$ has at most linear growth in space, locally uniformly in $t\ge 0$. 


We consider in this section the case of an autonomous volatility coefficient defined by 
\[
\forall\, (t,x)\!\in \R_+\times \R,\qquad\sigma(t,x) = \sigma(x).
\]
Note  that if  the solution $(X_t)_{t\ge 0}$ to~\eqref{eq:Volterrameanrevert}  with $X_0\!\in L^2(\P)$ is stationary~(\footnote{in the sense that the shifted processes $(X_{t+u})_{u\ge 0}$ and $(X_u)_{u\ge 0}$ have the same distribution as processes on the canonical space ${\cal C}(\R_+,\R)$.})  then $\mu(t)= \mu_0$ $d$-$a.e.$ and both  the mean   and  the variance of $X_t$ are  constant and finite as functions of $t$ as well as the expectations of any function of $X_t$  with at most quadratic growth. Typically such is the case of $x$, $x^2$ and $\sigma^2(x)$.

As a consequence, (see~\eqref{eq:L^p-supBound}) 
$$
\forall\, t\ge 0, \quad \E\, X_t = \frac{\mu_0}{\lambda}, \quad {\rm Var}(X_t)= v_0 \ge 0\quad \mbox{ and }\quad \bar \sigma_0^2 := \E\, \sigma^2(X_t) \ge 0.
$$
In fact as soon as $\mu(t) = \mu_0$ , Equation~\eqref{eq:Volterrameanrevert} writes 
The theorem below shows what are the consequences when only taking account  these three constraints

\begin{Theorem}[Time homogenous $\sigma$ and $\mu(t)=\mu_0$ is constant] \label{thm:atonome} Set $\mu(t)=\mu_0$ and $\sigma(t,x)= \sigma(x)$ in~\eqref{eq:Volterrameanrevert}. Assume $X_0\!\in L^2(\P)$ and let $(X_t)_{t\ge 0}$ be a solution  to~\eqref{eq:Volterrameanrevert}  satisfying 
\[
\forall\, t\ge 0, \quad \E\, X_t = \frac {\mu_0}{\lambda},\quad   {\rm Var}(X_t)=v_0\ge 0\quad \mbox{and} \quad   \E\, \sigma^2(X_t) = \bar \sigma_0^2\ge 0.
\]
%
\noindent Then, only two situations can occur:

\medskip 
\noindent $(i)$ If $\bar \sigma >0$, then the kernel $K$ is constant  with $K= \frac{2\lambda v_0}{\bar \s^2}$, so that~\eqref{eq:Volterrameanrevert} is a standard Brownian SDE and $(X_t)_{t\ge 0}$  is a (Markov) Brownian diffusion process with constant mean and variance.

\medskip 
 \noindent $(ii)$ If $\bar \sigma =0$, then $\sigma\big( \frac{\mu_0}{\lambda}\big)=0$ and $X_t =  \frac{\mu_0}{\lambda}$ $\P$-$a.s.$
\end{Theorem}
 
 \noindent {\bf Proof.} First note that as $\mu(t)= \mu_0 $ and $ R'_{\lambda}=-f_{\lambda} $, 
 Equation~\eqref{eq:Volterrameanrevert2} reads 
 \[
 X_t-\frac{\mu_0}{\lambda} = \Big(X_0-\frac{\mu_0}{\lambda}\Big) R_{\lambda}(t)+  \frac{1}{\lambda}\int_0^t f_{\lambda}(t-s)\sigma(X_s)dW_s.
 \]
 By  It\^o's isomorphism and Fubini's Theorem
\begin{align*}
{\rm Var}\Big(\int_0^t f_{\lambda}(t-s)\sigma(X_s)dW_s \Big) &= \E\Big(\int_0^t f_{\lambda}(t-s)\sigma(X_s)dW_s \Big)^2\\
& = \int_0^t f_{\lambda}(t-s)^2  \E\, \sigma(X_s)^2ds  =  \bar \sigma_0^2\int_0^t f_{\lambda}(t-s)^2 ds . 
\end{align*}
Then, it follows from the form~\eqref{eq:Volterrameanrevert2} of the scaled Volterra equation that 
\begin{equation}\label{eq:Station0}
\forall\, t\ge 0, \quad   v_0 = {\rm Var}(X_t)= v_0R^2_{\lambda}(t) + \frac { \bar \sigma^2}{\lambda^2 }\int_0^t f_{\lambda}^2(s)ds
\end{equation}
or, equivalently
\begin{equation}\label{eq:VolterraVar}
v_0\big(1-R^2_{\lambda}(t)\big) =  \frac { \bar \sigma^2}{\lambda^2 }\int_0^t f_{\lambda}^2(s)ds.
\end{equation} 
%
%
$(i)$ If $ \bar \sigma^2\neq 0$, differentiating this equality implies, since $R'_{\lambda} = - f_{\lambda}$ and  $R_{\lambda}$ and   $f_{\lambda}$ are  $>0$ on $(0,+\infty)$~--~which is a by-product of Assumption~$({\cal K})$, see~\eqref{eq:hypoRlambda})~--~that
\[
 \forall\, t>0, \quad  R'_{\lambda}(t)= -\kappa\, R_{\lambda}(t)  \quad \mbox{with $\kappa= \frac{2\lambda^2v_0}{\bar  \sigma^2}$}.
\]
Hence
\[
 \forall\, t\ge 0,\quad R_{\lambda}(t)= R_{\lambda}(0) e^{-\kappa t}= e^{-\kappa t}.
\]

Note  this also reads equivalently in terms of  in terms  of Laplace transform, $L_{R_{\lambda}}(t)= \frac{1}{\kappa+t}$, $t>0$.
Thus resolvent Equation~\eqref{eq:Resolvent2}, rewritten in terms  of Laplace transform, reads
\[
L_K(t)=  \frac{1}{\lambda}\Big( \frac{1}{tL_{R_{\lambda}}(t)}-1\Big)= \frac{\kappa}{\lambda}\frac{1}{t},\quad t>0
\]
which implies
\[
K(t) = \frac{\kappa}{\lambda}= \frac{2 \lambda v_0}{\bar  \sigma^2} ,\quad t>0.
\]
Hence  the kernel $K$ is necessary constant which yields the announced conclusion.

\medskip
\noindent $(ii)$ If $\bar \sigma^2=0$, then $v_0=0$ since $R_{\lambda}(t)<1$  for $t$ large enough. Consequently ${\rm Var}(X_t)=0$ for every $t\ge 0$. As $\E\, X_t = \E\, X_0=  \frac{\mu_0}{\lambda}$ owing to what was done in  section~\ref{sec:moyennecste}, it follows that  $X_t= \frac{\mu_0}{\lambda}$ $\P$-$a.s.$. But then $\sigma^2 \big( \frac{\mu_0}{\lambda} \big)= \bar \sigma^2 = 0$. \hfill$\Box$

\medskip
\noindent {\bf Examples.} $\rhd$ {\em The case $\sigma(x)=\sigma>0$ and $\mu(t)= \mu_0$.}  As $\bar \sigma =\sigma>0$, the existence of a stationary regime implies that the kernel $K$ is constant equal to $K(t)=\frac{2\lambda^2 v_0}{\sigma^2}$. Note that forcing $K=\mbox{\bf1}$ yields the expected constraint $v_0= \frac{\sigma^2}{2\lambda}$.
The process $X$ is then a standard  mean-reverting Ornstein-Uhlenbeck (Markov) process 
%
%
%
%
\noindent whose stationary regime does exist, is unique, with ${\cal N}\Big(\tfrac{\mu_0}{\lambda},v_0\Big)$ as  a $1$-marginal distribution (a.k.a. an invariant distribution in such a Markovian setting). This regime is of course obtained  by simply considering $X_0\sim {\cal N}\Big(\tfrac{\mu_0}{\lambda}, v_0\Big)$,  independent of $W$ so that $X$ is then a Gaussian process.

\medskip
\noindent $\rhd${\em Volterra $CIR$}.  This  is the class of processes of the   form
\[
X_t =X_0 +\int_0^t K(t-s)(\mu_0-\lambda X_s)ds +\int_0^t K(t-s) \vartheta \sqrt{X_s} dW_s, \quad X_0\ge 0
\]
with $\E\, X_0 = \frac{\mu_0}{\lambda}$, where $\vartheta, \lambda>0$ and $\mu_0> 0$. These processes are extensively used in rough volatility modeling with a kernel $K$ given by $K(u) =K_{H+\frac12}(u)=  u^{H-\frac 12}$, $H\!\in (0, \frac 12)$(see~\cite{JaissonR2016,ElEuchR2018,GatheralJR2018,GaJuRo2020, Jacquieretal2022}), provided the equation has at least  a non-negative weak solution  as established e.g.  in~\cite{JaissonR2016}. However, this specific  choice of $K$ plays no role here.
%
For this dynamics, 
$$
\bar \sigma_0^2=\bar \sigma(t) ^2 = \vartheta^2\E\, X_t =  \vartheta^2\frac{\mu_0}{\lambda}
$$
is constant whatever the choice of the kernel is  since $\E \, X_t$ is constant by assumption.

\smallskip
--  If $\bar \sigma_0^2>0$, then we are in the  case where the kernel $K$ is constant  equal to $\frac{2\lambda^2v_0}{\mu_0\vartheta^2}$ i.e. the regular $CIR$ model whose $1$-marginal invariant distribution is known to be the gamma distribution $\gamma\big(2\frac{\mu_0^2}{\lambda^2},\frac{\lambda v_0}{\mu_0} \big)$ with our notations. 

\smallskip
-- If $\bar \sigma_0^2=0$, it follows   that $\bar \sigma_0^2 = \frac{\mu_0}{\lambda}$ so that $\mu_0 = 0$. In turn~\eqref{eq:VolterraVar}  implies $v_0=0$ since $R_{\lambda}(t)< 1$ for $t$ large enough. If such is the case, then the trivial process $X_t= 0$, $t\ge 0$, is the only  stationary solution.

\section{Fake stationarity for  scaled Volterra equations}\label{sec:TowardFake}
In this section we investigate the case where, for every $t\ge 0$,   
\begin{equation}\label{eq:varsigma-sigma}
\mu(t)= \mu_0\quad \mbox{ and }\quad \sigma(t,x)=\varsigma(t) \sigma(x),\quad  \varsigma(t), \;\sigma(x)>0
\end{equation} 
\noindent where $\varsigma$ is a Bounded Borel function to be specified later. Then Equation~\eqref{eq:Volterrameanrevert}  writes
\begin{equation}\label{eq:Volterrameanrevert3}
 X_t =\frac{\mu_0}{\lambda} + \Big(X_0-\frac{\mu_0}{\lambda}\Big) R_{\lambda}(t)+  \frac{1}{\lambda}\int_0^t f_{\lambda}(t-s)\varsigma(s)\sigma(X_s)dW_s.
 \end{equation} 

Note that if the kernel $K$ satisfies~\eqref{eq:Kint} and~\eqref{eq:Kcont}, $\sigma$ is {\em Lipschitz continuous} and $\varsigma$ is a {\em bounded Borel function}, the above scaled Volterra equation  has a unique  ${(\cal F}^{X_0,W}_t)_{t\ge 0}$-adapted pathwise solution starting from $X_0\!\in L^2(\P)$ independent of $W$ (still owing to~\cite[Theorem~1.1]{JouPag22} once a reverse Wiener-Hopf transform has been performed, see Proposition~\ref{prop:Wienr-HopfX}). 

\medskip
\noindent {\bf Notational warning !} Also note that from now on we will still denote  $\bar\sigma^2(t) =\E\, \sigma^2(X_t)$.

\medskip
We saw that a {\em necessary} condition for the existence of a stationary regime is that $\E\, X_t$ and ${\rm Var}(X_t)$ are  constant  but also $\E\, \sigma^2(X_t)$, $t\ge 0$. 
The converse is clearly false as highlighted in the previous section.  However the existence of   regimes with constant means and variance remains interesting in practice to produce somewhat stabilized models. 

 This leads us to introduce a new object, the {\em stabilizer} and two notions of {\em fake stationary regimes} for a solution to the above Volterra equation~\eqref{eq:Volterrameanrevert3}.

\subsection{Stabilizer and fake stationary regimes}
\noindent $\blacktriangleright$  {\em Fake stationary regime of type~I}:  The mean $\E\, X_t$, the variance ${\rm Var}(X_t)$ 
and   $ \bar \sigma^2= \E\, \sigma^2(X_t)$ are all  
  constant as  functions of $t$.
  
 \medskip
Asking $ \E\, \sigma^2(X_t)$ to constant may appear as a more technical assumption that can be discussed. However it can also be seen as an  additional constraint which is one more step toward stationarity. 
We will see in Proposition~\ref{prop:equiv} further on that imposing that $\bar \sigma^2(t)$ is constant turns out to be superfluous as it can be obtained as a by-product of the constance of the first two quantities.

\medskip
\noindent $\blacktriangleright$  {\em Fake stationary regime of type~II}:   The solution $(X_t)_{t\ge 0}$ has the same marginal distribution i.e. $X_t\stackrel{d}{=} X_0$ for every $t\ge 0$.

\smallskip
At this stage, we can introduce the notion of {\em stabilizer}. Mimicking the  proof of Theorem~\ref{thm:atonome} where  the time-homogenous  $\sigma (x)$ is replaced by $\varsigma(t)\sigma(x)$  yields a {\em necessary condition} for a fake stationary regime of type~I to exist, namely  that
 the triplet    $\big(v_0, \bar  \sigma_0^2,\varsigma^2(t)\big)$ where $v_0 = {\rm Var}(X_0)$, $\bar \sigma_0^2 =\E\, \sigma(X_0)^2$  satisfies the equation
\begin{equation}\label{eq:VolterraVarbter}
\forall\, t>0, \quad c\lambda^2\big(1- R_{\lambda}^2(t)\big) = (f_{\lambda}^2*  \varsigma^2)_t
\end{equation}
 with $ \displaystyle c=\frac{v_0}{\bar \sigma_0^2}>0$.
%

\begin{Definition}[Corrector/Stabilizer]\label{def:corrrector}
We will  call {\em stabilizer} (or  {\em corrector} of  the scaled Volterra equation~\eqref{eq:Volterrameanrevert3}  bounded Bore function  $\varsigma= \varsigma_{\lambda, c}$ solution to the above equation~\eqref{eq:VolterraVarbter} (if any) . 
\end{Definition}

\noindent {\bf Remark.} $\bullet$ Note that~\eqref{eq:VolterraVarbter} has a solution $\varsigma_{\lambda,c}$  for some $c>0$ iff it has a solution $\varsigma_{\lambda,1}$ when $c=1$ and  $ \varsigma_{\lambda,c}= \sqrt{c} \,\varsigma_{\lambda,1}$.

\smallskip
\noindent $\bullet$ The existence of a stabilizer does no imply the existence of a fake stationary regime (type~I) for the Volterra equation with diffusion coefficient $\sigma(t,x)= \varsigma_{\lambda,c}(t)\sigma(x)$.

\medskip
This equation can be rewritten formally  in terms of Laplace transform
\[
c\lambda^2L_{1- R_{\lambda}^2} = L_{f_{\lambda}^2}L_{\varsigma^2}.
\]
In order to get rid of the Laplace transform of $1-R^2_{\lambda}$, we perform an integration by parts  
\[
L_{1-R^2}(t) = L_{1}(t)-L_{R_{\lambda}^2}(t)= \frac 1t -\Big(\frac{R^2_{\lambda}(0)}{t} -\frac 2t L_{R_{\lambda}f_{\lambda}}(t)\Big)= \frac 2t L_{R_{\lambda}f_{\lambda}}(t)
\]
so that the  Laplace counterpart  of~\eqref{eq:VolterraVarbter} finally reads

  \begin{equation}\label{eq:Laplacesigma}
\forall\, t>0, \quad  t\,L_{f^2_\lambda}(t).L_{\varsigma^2}(t)=  2\,c\lambda^2 L_{R_{\lambda}f_{\lambda}}(t).
\end{equation}

\medskip 
Note that, c being fixed, the solution $\varsigma_{\lambda,c}^2$ of~\eqref{eq:VolterraVarbter}, if any, is unique. Indeed, as $f_{\lambda}$ is a probability density owing to Assumption~$({\cal K})(ii)$ (see~\eqref{eq:hypoRlambda}), $L_{f^2_{\lambda}}>0$.  Then $L_{\varsigma^2}$ is uniquely determined by~\eqref{eq:Laplacesigma} which in turn {\it implies uniqueness} of $\varsigma^2$ (at least $dt$-a.e.). 

To carry on our investigations of the fake stationary regimes we make the following   existence assumption on Equation~\eqref{eq:VolterraVarbter}, in particular in terms of positivity: let $\lambda$, $c>0$,
\[
(E_{\lambda}) \quad  \mbox{Equation~\eqref{eq:VolterraVarbter} has a unique positive bounded  Borel solution $\varsigma_{\lambda,1}$ on $[0,+\infty)$ when $c=1$}.
\]

 Uniqueness is not a difficult point in practice (see e.g. the proof of the proposition below). The key condition in $(E_{\lambda})$ is  positiveness. 

 \medskip
 \noindent {\bf Remark.}  If $K=1$, then $R_{\lambda}(t)= e^{-\lambda t}$ and $f_{\lambda}(t)= \lambda  e^{-\lambda t}$ so that $\varsigma_{\lambda, c}= \sqrt{2\lambda c}$.
 
%
 
 \bigskip
 Numerical aspects  of this equation in the ``rough setting'', that is  for kernels $K= K_{\alpha}$, $\alpha = H+\frac 12$ with $H\!\in (0,\frac 12)$,  and for a specific family of diffusion coefficients $\sigma(x)= \sqrt{\kappa_0 + \kappa_1(x-\frac{\mu_0}{\lambda})+\kappa_2(x-\frac{\mu_0}{\lambda})^2}$, 
 are  investigated   in Section~\ref{sec:sigma2rough} where theoretical results are established like the  square integrability of $f_{\lambda}$ and  the existence of a continuous bounded stabilizer  $\varsigma_{\lambda,c}$, null at $0$ and positive on $(0, +\infty)$


\subsection{Fake stationary regimes of types~I and~II: first properties and examples}
Before exploring the main two situations for $\sigma$, let us come back to the notion of fake stationary regime of type~I.  The following proposition emphasizes that there is in fact an equivalence between having a constant variance and a constant value for $\E\, \sigma^2(X_t)$.

\begin{Proposition}[An equivalence] \label{prop:equiv}  Let $\lambda >0$, let $\mu_0\!\in \R$ and let $\s:\R\to \R
$ be a Lipschitz continuous function. Let  $X_0\!\in L^2(\Omega, {\cal A}, \P)$ be such that $\E\, X_0 = \frac{\mu_0}{\lambda}$ and  ${\rm Var}(X_0)=v_0\ge 0$.   Set $\bar \sigma_0^2 = \E\, \sigma^2(X_0)>0$ and set $c= \frac{v_0}{\bar \sigma_0^2}\!\in \R_+$. Assume the kernel     $K$ satisfies~\eqref{eq:Kint},~\eqref{eq:Kcont}, $f^2_{\lambda} $ has a finite Laplace transform on $(0, +\infty)$ and $(E_{\lambda})$ is in force.               

Then,  the unique strong solution starting from $X_0$ of the scaled Volterra equation~\eqref{eq:Volterrameanrevert3}  where 
$\varsigma_{\lambda,c}$ is solution to~\eqref{eq:VolterraVarbter} (has constant mean) and) satisfies the  following equivalence between the following two claims:

\medskip
\noindent $(i)$ $\forall\, t\ge 0$, ${\rm Var}(X_t )= {\rm Var} (X_0)=v_0$,

\medskip
\noindent $(ii)$ $\forall\, t\ge 0$, $\E\, \sigma^2(X_t )=  \E\, \sigma^2(X_0)= \bar \s_0^2$.
\end{Proposition} 

\noindent {\bf Proof.} We just have to prove the equivalence. 

\smallskip 
\noindent \fbox{$(ii)\Rightarrow (i)$} Then $\bar \sigma^2 =\bar \sigma^2_t = \E\, \sigma^2(X_t)$, $t\ge 0$, so that, having in mind that $c = \frac{v_0}{\bar \sigma^2_0}=    \frac{v_0}{\bar \sigma^2}$, it follows from Equation~\eqref{eq:VolterraVarbter} that
\begin{align*}
{\rm Var}(X_t ) &= {\rm Var}(X_0 )R^2_{\lambda}(t) + \frac{\bar \sigma^2}{\lambda^2} f^2_{\lambda}*\varsigma_{\lambda,c}^2 \\
&= {\rm Var}(X_0 )R^2_{\lambda}(t) + \frac{v_0}{c\lambda^2} f^2_{\lambda}*\varsigma_{\lambda,c}^2\\
& =  v_0 R^2_{\lambda}(t) +\big(1- R^2_{\lambda}(t) \big) v_0 = v_0.
\end{align*}
\noindent \fbox{$(i) \Rightarrow (ii)$} Assume  $  {\rm Var}(X_t) = {\rm Var}(X_0) = v_0$ for every  $t\ge 0$. 

If $v_0=0$, then $X_t = \frac{\mu_0}{\lambda}$ $a.s.$ for every $t\ge 0$. Consequently $\E\, \sigma^2(X_t )=  \sigma^2\big(\frac{\mu_0}{\lambda}\big)$, $t\ge 0$. 

If $v_0>0$, $c\!\in (0, +\infty)$ since  $\bar \sigma^2_0 = \E\,\sigma^2(X_0)>0$. 
Then, one checks using~\eqref{eq:Volterrameanrevert2} that the function  
\[
y(t) = \varsigma^2_{\lambda,c}(t)\Big( \frac{\E\,\sigma^2(X_t)}{\bar \sigma^2_0}-1\Big)
\]
is solution to $f^2_{\lambda} * y=0$ and satisfies $ y(0) = 0$ since $\varsigma_{\lambda,c}$ is bounded by assumption and $t\mapsto \E\, \sigma^2(X_t)$ is continuous. Moreover by the linear growth assumption on $\sigma$, we know that $\E\, \sigma^2(X_t)$ is bounded since $\E\, X_t^2$ is.  Consequently $y$ has a Laplace transform as well as its positive and negative parts $y^\pm$. One has $L_{f^2_{\lambda}} \cdot L_{y^+} = L_{f^2_{\lambda}} \cdot L_{y^-} $. If $L_{f^2_{\lambda}}$ is not identically $0$, then it is positive on $(0, +\infty)$ so that $L_{y^+}= L_{y^-}$ which implies $y^+=y^-$, hence $y=0$. \hfill $\Box$

 \medskip
 This result does not presume of the existence of a fake stationary regime of type~I, but that  to establish its existence, it suffices to check that both mean and variance are constant.

\subsubsection{The case $\sigma(x)= \sigma$} Then under the assumptions of Proposition~\ref{prop:equiv}  one easily checks that  $\bar \sigma_0^2 = \sigma^2$ and that  the solution $(X_t)_{t\ge 0}$  to~\eqref{eq:Volterrameanrevert3} with $c= \frac{v_0}{\sigma^2}$ has  constant mean $\frac{\mu_0}{\lambda}$ and variance $v_0$. As a consequence, if $X_0\sim {\cal N}\big(\frac{\mu_0}{\lambda},v_0\big)$, then $X_t \sim X_0$  for every $t\ge 0$ i.e. $X$  has a  fake stationary regime of type~II since the process becomes clearly Gaussian. All this follows from the identity: 
\[
\forall\, s,\, t\ge 0, \quad {\rm Cov}(X_t,X_s) = v_0\Big(R_{\lambda}(t)R_{\lambda}(s) + \frac {1}{\lambda^2c} \int_0^{t \wedge s} f_{\lambda}(t-u)f_{\lambda}(s-u)\varsigma^2_{\lambda,c}(u)du\Big).
\]

\subsubsection{Trinomial squared volatility}\label{sec:FakeI-II}

In this section we specify a family of stabilized scaled Volterra equations~\eqref{eq:Volterrameanrevert3}  where 
\begin{equation} \label{eq:sigmafakeI&II}
\sigma(x) = \sqrt{ \kappa_0 +\kappa_1\,(x-\tfrac{\mu_0}{\lambda})+\kappa_2\,(x-\tfrac{\mu_0}{\lambda})^2}\quad \mbox{ with }\quad  \kappa_i\ge 0,\;i=1,2, \;\kappa^2_1 \le 4\kappa_0\kappa_2,
\end{equation}
for which we will establish the existence of a fake stationary regime (of type~I in general). Note that  $\kappa_1=0$ whenever $\kappa_0$ or $\kappa_2=0$ and that
\[
\kappa_2 = [\sigma]^2_{\rm Lip} \quad \mbox{ and }\quad \sigma\big( \tfrac{\mu_0}{\lambda}\big) = \kappa_0.
\] 
This  is this type of vol-vol term that appears in the quadratic rough volatility model introduced in~\cite{GaJuRo2020} to solve the problem of the joint S\&P 500/VIX smile calibration problem.

\begin{Proposition}[Fake stationary regimes  (types~I and II) and first asymptotics]\label{prop:quasiZumbach} Consider the stabilized Volterra Equation~\eqref{eq:Volterrameanrevert3} with $\lambda>0$, $\mu_0\!\in \R$,  where $\sigma$ is given by~\eqref{eq:sigmafakeI&II} and   $\varsigma = \varsigma_{\lambda,c}$,  assumed to be the unique continuous solution to Equation~\eqref{eq:VolterraVarbter} for some $c \in (0, \frac{1}{\kappa_2})$ (so that condition $(E_{\lambda})$ is satisfied).
Then, let $X_0\!\in L^2(\P)$ with 
\[
\E\, X_0 = \frac{\mu_0}{\lambda} \quad \mbox{ and }\quad v_0= {\rm Var}(X_0) = \frac{c\kappa_0}{1-c\kappa_2}.
\]

\noindent $(a)$ {\em Case $\kappa_2>0$}. The solution $ (X_t)_{t\ge 0}$  to~\eqref{eq:Volterrameanrevert3} has a fake stationary regime  of type~I  i.e. 
\[
\forall\, t\ge 0, \qquad \E\, X_t = \frac{\mu_0}{\lambda},  \quad{\rm Var}(X_t) = v_0=\frac{c\kappa_0}{1-c\kappa_2} \mbox{ and }\quad \E\, \sigma^2(X_t) =  \bar \sigma_0^2 =\frac{\kappa_0}{1-c\kappa_2}.
\]
Moreover, for any starting value $X_0\!\in L^2(\P)$, 
\[
 \E\, X_t \to  \frac{\mu_0}{\lambda} \quad \mbox{ and }\quad {\rm Var}(X_t) \to  \frac{c\kappa_0}{1-c\kappa_2}\quad \mbox{ as }t\to +\infty.
\]

\noindent $(b)$ {\em Case $\kappa_2=0$}. Then 
$\kappa_1=0$, $\sigma(x)= \sqrt{\kappa_0}$ is constant (the choice of $c>0$ becomes free) and if $X_0\sim \nu:={\mathcal N}\big( \frac{\mu_0}{\lambda}, \frac{c\kappa_0}{1-c\kappa_2}\big)$, then $ (X_t)_{t\ge 0}$  is a Gaussian process  with a  fake stationary regime of type~II  with $\nu$   as  $1$-marginal distribution 
 \end{Proposition}

\noindent {\bf Practitioner's corner}. $\bullet$ The above proposition covers at least partially  the rough volatility dynamics  recently introduced in~\cite{GaJuRo2020} (the volatility process being defined as $V_t = \sigma(X_t)$ with our notations). In  this model, the traded asset and its volatility are driven by t{\em he same Brownian motion} and its purpose is to propose a joint calibration of the {\em VIX} and the {\em SP500} in order to take into account  the so-called Zumbach effect which connects the  evolution of the asset (here an index) and its volatility.

\smallskip
\noindent $\bullet$ Note that in practice, when $\kappa_1>0$, if we fix the value of $v_0$, then $c = \frac{v_0}{\kappa_0+ v_0\kappa_2}$ so that, $\sigma$ being $\sqrt{\kappa_2}$-Lipschitz continuous, one has $\kappa_2[\sigma]^2_{\rm Lip}= \frac{v_0\kappa_2}{\kappa_0+ v_0\kappa_2} <1$ which ensures the $L^2$-confluence of the paths of the solution (see Proposition~\ref{prop:confluence} further on). 


\bigskip
\noindent {\bf Remark.}  If we consider a dynamics more inspired by the original CIR model in which $\mu(t) = \mu_0>0$, $\sigma(t,x)= \varsigma_{\lambda,c}(t)\sigma(x)$ with $\sigma(x) :=\sqrt{\kappa_0+\kappa_1(x-\frac{\mu_0}{\lambda})}$, $\kappa_1>0$, $\kappa_0>\kappa_1\frac{\mu_0}{\lambda}$, the resulting Volterra equation~\eqref{eq:Volterrameanrevert} may have (at least) a non-negative weak solution e.g. as the $C$-weak limit of a system of Hawkes processes  (as it has been done in~\cite{JaissonR2016} without the presence  of the corrector $\varsigma_{\lambda,c}$, see also~\cite{Jaber2017},~\cite{alfonsi2023}). If such a non-negative weak solution on the whole non-negative  real line does exists starting from $X_0$ with mean $\frac{\mu_0}{\lambda}$    then
\[
{\rm Var}(X_t) = {\rm Var}(X_0)R^2_{\lambda}(t) + \kappa_0\big(1-R^2_{\lambda}(t)\big).
\] 
A fake stationary regime of type~I then should have  constant mean $\E\, X_t = \frac{\mu_0}{\lambda}$ and a variance given by  ${\rm Var}(X_t) = c \kappa_0$ respectively and one easily checks that such is the case. Note that if $\kappa_0=0$~--~like  in the Volterra CIR like model~--~one retrieves the degenerate situation where $X_t= \frac{\mu_0}{\lambda}$ $\P$-$a.s.$

\bigskip
\noindent {\bf Proof.}  $(a)$ Note that, if $ \kappa_2>0$,  $[\sigma]_{\rm Lip} = \sqrt{\kappa_2}$ since 
\[
\sigma(x) = \sqrt{ \kappa_2\big(x-\tfrac{\mu_0}{\lambda}-\frac{\kappa_1}{2\kappa_2}\big)^2 +\kappa_0 -\frac{\kappa_1^2}{4\kappa_2}}
\]
so that the  condition $0<c<1/[\sigma]^2_{\rm Lip}$ is satisfied. We know that 
$$
\E\, X_t = \frac{\mu_0}{\lambda} R_{\lambda}^2(t) + \frac{\mu_0}{\lambda}(1-R^2_{\lambda}(t)) = \frac{\mu_0}{\lambda}
$$  
for every $t\ge 0$. As for the variance, we have
\begin{align}
\nonumber {\rm Var}(X_t) &= {\rm Var}(X_0) R^2_{\lambda}(t)  + \frac{1}{\lambda^2}f^2_{\lambda}* \big( \varsigma^2\,\E\, \sigma^2(X_{\cdot})\big)_t\\
\label{eq:VarStatio}& =  {\rm Var}(X_0) R^2_{\lambda}(t)  + \frac{\kappa_0 }{\lambda^2}(f^2_{\lambda}*\varsigma^2)_t +\kappa_1\times 0+  \frac{ \kappa_2}{\lambda^2}f^2_{\lambda}*\big(\varsigma^2 \, {\rm Var}(X_{\cdot})\big)_t.
\end{align}
If we assume that ${\rm Var}(X_t)$ is constant  i.e. ${\rm Var}(X_t)=v_0$ for every $t\ge 0$ and take advantage of the identity~\eqref{eq:VolterraVarbter} satisfied by $\varsigma=\varsigma_{\lambda,c}$, the above equation reads
\[
v_0(1-R^2_{\lambda}(t))= (c\kappa_0+c\kappa_2v_0)(1-R^2_{\lambda}(t)) , \quad t\ge 0
\]
i.e. 
\[
v_0 = \frac{c\kappa_0}{1-c\kappa_2}>0
\] 
since $R_{\lambda}(t)>1$ for $t$ large enough ,which is clearly solution to the equation. 

Conversely one checks that this constant value for the variance solves the above equation. Let us prove that it is the only one. By the linearity of Equation~\eqref{eq:VarStatio}, it suffices to show  that  the equation in $x\!\in {\cal C}(\R_+,\R)$
\[
x(t) =  \frac{\kappa_2}{\lambda^2} \big(f^2_{\lambda}*(\varsigma^2.\, x)\big)_t ,\quad x(0)=0
\]
only has  the null function as solution. If $x$ solves the above equation, then 
\[
|x(t)|\le  \frac{\kappa_2}{\lambda^2} (f^2_{\lambda}*\varsigma^2)_t  \sup_{0\le s \le t} |x(s)| \, = c\kappa_2    \sup_{0\le s \le t} |x(s)| .
\]
If $x\equiv \!\!\!\!\! / \,0$, there exist $\varepsilon>0$ such that $\tau_{\varepsilon} = \inf\{t: |x(t)|>\varepsilon \} <+\infty$. By continuity of $x$ it is clear that $\tau_{\varepsilon}  >0$ and $ |x(\tau_{\varepsilon})|= \sup_{0:\le s \le t} |x(s)| = \varepsilon$ which is impossible since $\kappa_2c<1$. Consequently $x\equiv 0$. 

Hence $(X_t)_{t\ge 0}$ is  a fake stationary regime of type~II with the above mean and variance. The last claim is a straightforward  consequence of Proposition~\ref{prop:confluence}.

\smallskip
\noindent $(b)$ is obvious  once noted that $(X_t)_{t\ge 0}$ is a Gaussian process (and $[\sigma]_{\rm Lip}=0$).
 \hfill$\Box$

\section{Long run behaviour : asymptotics and confluence}
In this section we try to answer the question; ``How does the process $(X_t)_{t\ge}$ behave when $t$ goes to infinity and possibly starting from any $X_0\!\in L^p(\P)$, $p=2$ or $p>2$~?''

\smallskip
We first recall a result on the constant of an $L^p$-Burkholder-Davis-Gundy (BDG) inequality. 
\begin{Lemma}[Best constant in  a BDG  inequality (see Remark~2 in~\cite{CarlenK1991})] Let $M$ be a continuous local martingale null at $t=0$. Then, for every $p\ge 1$ 
\[
\|M_t\|_p \le 2\sqrt{p} \,\| \langle M\rangle^{\frac 12}\|_p.
\] 
\end{Lemma}

\begin{Proposition}[Moment control] \label{prop:Momentctrl} Let $\lambda,c>0$, let $\mu_0\!\in \R$, and let $\sigma:\R\to \R$ be a Lipschitz continuous function. Assume $(E_{\lambda})$ is in force  and $\varsigma_{\lambda,c}$ is continuous. Assume $f_{\lambda}\!\in {\cal L}^2({\rm Leb}_1)$. Let $(X_t)_{t\ge 0}$ be a solution to the scaled Volterra equation in its form~\eqref{eq:Volterrameanrevert3} starting from  any $X_0\!\in L^2(\P)$.

\smallskip
\noindent $(a)$ {\em First two  moments}.  Assume $c\!\in \big(0,\frac{1}{[\sigma]^2_{\rm Lip}}\big)$.
Then,  one has 
\[
\Big| \E \,X_t-\frac{\mu_0}{\lambda}\Big| \le R_{\lambda}(t) \Big| \E\,X_0-\frac{\mu_0}{\lambda}\Big| , \quad t\ge 0
\]
and 
\[
 \sup_{t\ge 0}\Big\|X_t-\frac{\mu_0}{\lambda}\Big\|_2 \le   \left[\frac{\sqrt{c}}{1-[\sigma]_{\rm Lip}\sqrt{c}}\big|\sigma(\tfrac{\mu_0}{\lambda})\big|\right]\vee   \Big\|X_0-\frac{\mu_0}{\lambda}\Big\|_2< +\infty.
\]
$(b)$ {\em $L^p$-moments, $p\!\in (2,+\infty)$.} 
If furthermore $X_0\!\in L^p(\P)$ and $c\!\in \big(0, \frac{1}{4p[\sigma]^2_{\rm Lip}}\big)$ 
 then 
\[
 \sup_{t\ge 0} \Big\|X_t-\frac{\mu_0}{\lambda}\Big\|_p\le \inf_{\epsilon\in (0, \frac{1}{4cp[\sigma]_{\rm Lip}^2}-1)} \bigg( \left[\frac{2\sqrt{pc(1+\epsilon)}}{1-2[\sigma]_{\rm Lip}\sqrt{pc(1+\epsilon)}}\big|\sigma(\tfrac{\mu_0}{\lambda})\big|\right]\vee \left[ (1+1/\epsilon)^{\frac12} \Big\|X_0-\frac{\mu_0}{\lambda}\Big\|_p\right]\bigg)< +\infty.
\]
\end{Proposition}

\noindent{\bf Remark.} If $|\sigma(x)|^2 \le \kappa_0 +\kappa_2(x-\frac{\mu_0}{\lambda})^2$ is a Borel function, the above results remain valid by simply replacing $[\sigma]^2_{\rm Lip}$ by $\kappa_2$ if~\eqref{eq:Volterrameanrevert}  has a solution.

\medskip 
\noindent {\bf Proof.} $(a)$ Set $\rho := c[\sigma]^2_{\rm Lip}\!\in (0,1)$ and let $\eta \in (0, 1-\rho)$ be a free parameter such that $\rho+\eta \!\in (0,1)$. 
One has for every $x\!\in \R$,
\begin{equation}\label{eq:quadrabound}
\sigma^2(x)   \le \big(\big|\sigma(\tfrac{\mu_0}{\lambda}) \big|+ [\sigma]_{\rm Lip}\big |x -\tfrac{\mu_0}{\lambda}|\big)^2\le\kappa_0+\kappa_2 \big |x -\tfrac{\mu_0}{\lambda}|^2
\end{equation}
with real constant $\kappa_i$, $i=0,2$ depending on $\eta$ and reading
$$
\kappa_0=  \kappa_0(\eta)=(1+\rho\eta^{-1}) |\sigma(\tfrac{\mu_0}{\lambda})|^2 \quad \mbox{ and} \quad \kappa_2=    \kappa_2(\eta)=[\sigma]_{\rm Lip}^2(1+\eta\rho^{-1}) 
 $$
so that $c\kappa_2= \rho+\eta<1$.
Using that $f^2_{\lambda}*\varsigma^2 = c \lambda^2(1-R^2_{\lambda})$ elementary computations show that for every $t\ge 0$
\begin{equation}\label{eq:majormomentL2}
\E\, \Big( X_t-\frac{\mu_0}{\lambda}\Big)^2 \le \E\, \Big( X_0-\frac{\mu_0}{\lambda}\Big)^2 R^2_{\lambda}(t) + \kappa_1\big( 1-R^2_{\lambda}(t) \big) + \frac{\kappa_2}{\lambda^2} \int_0^t f^2_{\lambda}(t-s) \varsigma^2(s) \E\, \Big( X_s-\frac{\mu_0}{\lambda}\Big)^2 ds
\end{equation}\label{eq:variance bound}
Now let $A> \bar A_{\eta}:=\frac{\kappa_0c}{1-\kappa_2c}\vee  \E\, \Big( X_0-\frac{\mu_0}{\lambda}\Big)^2$, $\delta>0$ and 
\[
t_{\delta} = \inf\Big\{ t:  \E\, \Big( X_t-\frac{\mu_0}{\lambda}\Big)^2 \ge A + \delta\Big\}.
\]
As $t\mapsto \E\, \Big( X_t-\frac{\mu_0}{\lambda}\Big)^2$  is continuous  and $A>\E\, \Big( X_0-\frac{\mu_0}{\lambda}\Big)^2$ it follows from the above inequality  and the identity satisfied by $\varsigma$ that, if $t_{\delta}<+\infty$,  
\[
A+\delta = \E\, \Big( X_{t_\delta}-\frac{\mu_0}{\lambda}\Big)^2< A R^2_{\lambda}(t_{\delta}) + \big(\kappa_0 c +\kappa_2c (A+\delta)\big)\big(1-R^2_{\lambda}(t_{\delta})\big).
\]
Now, we have $\kappa_0 c +\kappa_2c A<A$ by construction of $A$, hence
\[
A+\delta = \E\, \Big( X_{t_{\delta}}-\frac{\mu_0}{\lambda}\Big)^2 < A R^2_{\lambda}(t_{\delta}) +A(1-R^2_{\lambda}(t_{\delta})) + \kappa_2c \delta\big((1-R^2_{\lambda}(t_{\delta})\big) < A + \delta.
\]
which yields a contradiction. Consequently, $t_{\delta}= +\infty$ which implies that $\E\, \Big( X_t-\frac{\mu_0}{\lambda}\Big)^2 \le A+\delta$ for every $t\ge 0$. Letting $\delta \to 0$ and 
$A\to \bar A_{\eta}$ successively, yields 
\[
 \sup_{t\ge 0}\E\, \Big( X_t-\frac{\mu_0}{\lambda}\Big)^2\le \bar A_{\eta}=  c |\sigma(\tfrac{\mu_0}{\lambda})|^2\frac{\eta+\rho}{\eta(1-\rho-\eta)}.
\]
Then one checks that $\eta \mapsto \bar A_{\eta}$ is minimal over $(0,1-\rho)$ at $\eta = \sqrt{\rho}-\rho$ which finally yields the announced result.

\smallskip
\noindent$(b)$ Let $p\ge 2$. Set  $\rho_p  := c\,C^{2}_p[\sigma]_{\rm Lip}^2<1$. Using  successively  {\em BDG} inequality to the (a priori) local martingale $M_s = \int_0^s f_{\lambda}(t-u)\zeta(u)\sigma(X_u)dW_u$, $0\le s\le t$, (see~\cite[Proposition~4.3]{RevuzYor}) and the generalized Minkowski inequality, we get
\begin{align*}
\Big\|  X_t-\frac{\mu_0}{\lambda}\Big\|_p& \le \Big\|  X_0-\frac{\mu_0}{\lambda}\Big\|_pR_{\lambda}(t) + \frac{C_p}{\lambda}\Big\|\int_0^t f^2_{\lambda}(t-s)\varsigma^2(s) \sigma(X_s)^2ds \Big\|_{\frac p2}^{\frac 12}\\
& \le  \Big\|  X_0-\frac{\mu_0}{\lambda}\Big\|_pR_{\lambda}(t) + \frac{C_p}{\lambda}\Big( \int_0^t f^2_{\lambda}(t-s)\varsigma^2(s) \|\sigma(X_s)\|^2_{p}\Big)^{\frac 12}
\end{align*}
 (with an equality when $p=2$ in the first line). 
  Now let $\epsilon\!\in (0,1/\rho_p-1)$. It follows from the elementary inequality $(a+b)^2 \le (1+1/\epsilon) a^2+(1+\epsilon )b^2$ that
 \[
 \Big\|  X_t-\frac{\mu_0}{\lambda}\Big\|_p^2 \le \Big\|  X_0-\frac{\mu_0}{\lambda}\Big\|^2_pR^2_{\lambda}(t)(1+1/\epsilon) + \frac{C^2_p}{\lambda^2}(1+\epsilon) \int_0^t f^2_{\lambda}(t-s)\varsigma^2(s) \|\sigma(X_s)^2\|_{\frac p2}ds.
 \]
 
 Let $\tilde\rho_p = \rho_p(1+\epsilon) \!\in (0,1)$ and for $\eta \!\in (0, 1-\tilde \rho_p)$, set now $\kappa_i$ depending on $\eta$ and reading
$$
\kappa_0=  \kappa_1(\eta) =(1+\tilde\rho_p\eta^{-1}) |\sigma(\tfrac{\mu_0}{\lambda})|^2 \quad \mbox{ and} \quad \kappa_2=   \kappa_2(\eta) =[\sigma]_{\rm Lip}^2(1+\eta\tilde\rho_p^{-1}) 
 $$
so that $c\,C^2_p(1+\epsilon)\kappa_2=\tilde\rho_p+\eta<1$.
  As $\frac p2\ge 1$, it follows from~\eqref{eq:quadrabound} that 
 \[
 \|\sigma(X_s)\|^2_{p} = \|\sigma(X_s)^2\|_{\frac p2} \le\kappa_0+\kappa_2\Big\| X_s-\frac{\mu_0}{\lambda}\Big\|^2_p
 \]
 which entails, combined with the identity  $f^2_{\lambda}*\varsigma^2 = c \lambda^2 (1-R^2_{\lambda})$,   that, for every $t\ge 0$, 
\begin{align}\nonumber 
 \Big\|  X_t-\frac{\mu_0}{\lambda}\Big\|_p^2 &\le  \Big\|  X_0-\frac{\mu_0}{\lambda}\Big\|_p^2 R^2_{\lambda}(t) (1+1/\epsilon)+C_p^2(1+\epsilon)\Big( \kappa_1c\big( 1-R^2_{\lambda}(t) \big) \\
\label{eq:variancebound2}
 & \qquad + \frac{\kappa_2}{\lambda^2} \int_0^t f^2_{\lambda}(t-s) \varsigma^2(s)\Big\| X_s-\frac{\mu_0}{\lambda}\Big\|^2_p
ds\Big).
\end{align}
Now let $A> \bar A_{\eta,\epsilon}:=\frac{\kappa_0c\,C_p^2(1+\epsilon)}{1-\kappa_2c\,C_p^2(1+\epsilon)}\vee  \Big[ (1+1/\epsilon)\Big\|  X_0-\frac{\mu_0}{\lambda}\Big\|_p^2\Big]$, $\delta>0$ and 
\[
t_{\delta} = \inf\Big\{ t: \Big\|  X_t-\frac{\mu_0}{\lambda}\Big\|_p^2 \ge A + \delta\Big\}.
\]
If $t_{\delta}<+\infty$, it follows on the one hand from the continuity  of $t\mapsto\Big\|  X_t-\frac{\mu_0}{\lambda}\Big\|_p^2$   that    $A+\delta =\Big\|  X_{t_\delta}-\frac{\mu_0}{\lambda}\Big\|_p^2$ and, on the other hand,   from  Equation~\eqref{eq:VolterraVarbter} satisfied by $\varsigma$, that
\[
\int_0^t f^2_{\lambda}(t-s) \varsigma^2(s)\Big\| X_s-\frac{\mu_0}{\lambda}\Big\|^2_pds\le A(1-R^2_{\lambda}(t)). 
\]
Moreover as $A>\Big\|  X_0-\frac{\mu_0}{\lambda}\Big\|_p^2(1+1/\epsilon)$, we deduce from~\eqref{eq:variancebound2}
\[
A+\delta =\Big\|  X_{t_\delta}-\frac{\mu_0}{\lambda}\Big\|_p^2< A R^2_{\lambda}(t_{\delta})+C_p^2(1+\epsilon)\big(\kappa_0 c +\kappa_2c (A+\delta)\big)\big(1-R^2_{\lambda}(t_{\delta})\big).
\]
Now, we have $C_p^2(1+\epsilon)c(\kappa_0   +\kappa_2 A)<A$ by definition of $A$. Hence
\[
A+\delta =\Big\|  X_{t_\delta}-\frac{\mu_0}{\lambda}\Big\|_p^2 < A R^2_{\lambda}(t_{\delta}) +  A(1-R^2_{\lambda}(t)) +   C_p^2(1+\epsilon)c\kappa_2 \delta\big((1-R^2_{\lambda}(t_{\delta})\big) < A + \delta.
\]
since $ C_p^2(1+\epsilon)c\kappa_2<1$. This  yields a contradiction. Consequently, $t_{\delta}= +\infty$ which implies that $\Big\|  X_{t_\delta}-\frac{\mu_0}{\lambda}\Big\|_p^2  \le A+\delta$ for every $t\ge 0$. Letting $\delta \to 0$ and 
$A\to \bar A_{\eta,\epsilon}$ successively, yields 
\[
 \sup_{t\ge 0}\Big\|  X_t-\frac{\mu_0}{\lambda}\Big\|^2_p\le \bar A_{\eta,\epsilon}<+\infty .
\]
Then one checks that $\eta \mapsto \bar A_{\eta,\epsilon}$ is minimal over $(0,1-\tilde \rho_p)$ at $\eta = \sqrt{\tilde \rho_p}-\tilde \rho_p$ which finally yields the announced result. \hfill$\Box$

\subsection{$L^2$-confluence}
The following result can be compared to the confluence property satisfied by he mean-reverting Ornstein-Uhlenbeck process (without 
rate of convergence however).

\begin{Proposition}[$L^2$-confluence]\label{prop:confluence}  Assume that all the above assumptions of Proposition~\ref{prop:Momentctrl} $(a)$ are in force. Assume furthermore that $f_{\lambda}\!\in L_{\R_+}^2({\rm Leb}_1)$.


\noindent $(a)$ There exists a continuous  non negative function $\varphi_{\infty}: \R_+\to [0,1]$ such that $\varphi_{\infty}(0)=1$, $\displaystyle \lim_{t\to +\infty} \varphi_{\infty}(t) =0$, only depending on $\lambda$, $c$ and the kernel $K$, such that  for any $X_0$,$X'_0\!\in L^2({\mathbb P})$,  the solutions to Volterra equation~\eqref{eq:Volterrameanrevert} denoted $(X_t)_{t\ge 0}$ and $(X'_t)_{t\ge 0}$ starting from $X_0$ and $X'_0$ respectively satisfy
\[
\forall\, t\ge 0, \quad  \E\, |X_t-X'_t|^2   \le \varphi_{\infty}(t)\, \E\,|X_0-X'_0|^2 .
\] 
\noindent $(b)$ In particular if $X$ has a fake stationary regime of type~I,   
\[
\forall\, t\ge 0, \quad  {\cal W}_2([X'_t], [X_t])  \le \varphi_{\infty}(t)^{1/2}\, {\cal W}_2([X'_0], [X_0])
\] 
so that $\E\, X'_t \to \frac{\mu_0}{\lambda}$,  ${\rm Var}(X'_t ) \to v_0$  and, more generally, finite dimensional ${\cal W}_2$-convergence.

\smallskip 
\noindent $(c)$  In case if $X$ has a fake stationary regime of type~II, its marginal distribution   is unique.

\smallskip 
\noindent $(d)$ In case $\sigma\big(\frac{\mu_0}{\lambda}\big)=0$ then $\big(X'_t = \frac{\mu_0}{\lambda}\big)_{t\ge0}$,  is a trivial stationary (degenerate) regime and, for every $X_0\!\in L^2({\mathbb P})$, the solution  $(X_t)_{t\ge 0}$  of~\eqref{eq:Volterrameanrevert} starting from $X_0$ satisfies
\[
\forall\, t\ge 0, \quad  \E\,\Big |X_t-\frac{\mu_0}{\lambda}\Big|^2   \le \varphi_{\infty}(t)\, \E\,\Big|X_0-\frac{\mu_0}{\lambda}\Big|^2.
\]
\end{Proposition}

\noindent {\bf Proof}. Set $\Delta _t = X_t-X'_t \!\in  L^2(\P)$ for every $ t\ge 0$.

\smallskip
\noindent $(b)$, $(c)$ and $(d)$ are all straightforward consequences of $(a)$.

\smallskip
\noindent $(a)$ We still denote $\rho = c[\sigma]_{\rm Lip}^2\!\in (0,1)$.   It follows from the reduced form~\eqref{eq:Volterrameanrevert2}, It\^o's isometry and the Lipschitz property of $\sigma$ that
\begin{align}
\label{eq:ineqDeltat} \E\, \Delta_t^2&  \le  R_{\lambda}^2(t)  \E\, \Delta_0^2+ \frac{[\sigma]_{\rm Lip}^2}{\lambda^2}\int_0^t f^2_{\lambda}(t-s)\varsigma^2(s) \E\, \Delta^2_s ds.
\end{align}

Let 
$\bar \delta_t = \| \Delta_t\|_2$ to alleviate notations.   One checks that, under our assumptions, the function $t\mapsto \bar \delta_t$ is continuous (see~\cite{JouPag22}).

Let $\eta>0$ such that $\rho(1+\eta)^2 <1$. We define $\tau_{\eta} :=\inf\{t: \bar \delta_t> (1+\eta)\bar \delta_0\}$. If 
$\tau_{\eta} <+\infty$, then  $\bar \delta_s \le  (1+\eta)\bar \delta_0$ for every $s\!\in (0, \tau_{\eta})$ and by continuity $\bar \delta_{\tau_{\eta}}^2= (1+\eta)^2\bar \delta_0^2$. Inserting this in the above inequality at time $\tau_{\eta}$ yields
\[
\bar \delta_{\tau_{\eta}}^2 \le \bar \delta_0^2 \big[R_{\lambda}^2(\tau_{\eta})+ (1- R_{\lambda}^2(\tau_{\eta}))\rho(1+\eta)^2\big] < \rho(1+\eta)^2\bar \delta_0^2
\]
which yields a contradiction so that $\bar \delta_s\le (1+\eta)\bar \delta_0$ for every  $s>0$. This  holds for every $\eta>0$  which in turn implies that $\bar \delta_t \le \bar \delta_0$ for every $t\ge 0$. Inserting this again in~\eqref{eq:ineqDeltat} implies that, for every $t>0$, 
\[
\delta^2_t \le \delta^2_0\varphi_1(t) \quad \mbox{ with }\quad \varphi_1(t):=R_{\lambda}^2(t)+(1-R_{\lambda}^2(t))\rho.
\]
Note that $\varphi_1(t) = \rho + R^2_{\lambda}(t)(1-\rho)$ satisfies 
\[
\varphi_1(0):=1, \; \varphi_1(t) \!\in (0,1), \, t>0 \mbox{ and }\; \varphi_1 \mbox{ is non-decreasing and continuous.}
\]
Inserting this upper-bounds of $\delta^2_t $ into~\eqref{eq:ineqDeltat} straightforwardly yields 
\[
\delta^2_t \le \bar \delta_0\varphi_2(t) \quad \mbox{ with }\quad \varphi_2(t):=R_{\lambda}^2(t)+ \rho\int_0^t f^2_{\lambda}(t-s)\varsigma^2(s) \varphi_1(s) \frac{ds}{\lambda^2 c}.
\]
One checks using  identity~\eqref{eq:VolterraVarbter}   satisfied by $\varsigma^2$ and the definition of $\varphi_1$ that 
\[
\varphi_2(t):=\varphi_1(t)- \rho\int_0^t f^2_{\lambda}(t-s)\varsigma^2(s)\big( 1- \varphi_1(s)\big) \frac{ds}{\lambda^2 c}
\]
so that $0\le \varphi_2< \varphi_1<1$ on $(0, +\infty)$. By induction one shows that 
\[
\delta^2_t \le \delta^2_0\varphi_k(t)
\]
with 
\begin{align*}
\varphi_{k}(t) &= R_{\lambda}^2(t)+ \rho\int_0^t f^2_{\lambda}(t-s)\varsigma^2(s) \varphi_{k-1}(s) \frac{ds}{\lambda^2 c}\\		    
 		     & =  \varphi_1(t)- \rho\int_0^t f^2_{\lambda}(t-s)\varsigma^2(s)\big( 1- \varphi_{k-1}(s)\big) \frac{ds}{\lambda^2 c}
\end{align*}
 where we used again~\eqref{eq:VolterraVarbter}  satisfied by $\varsigma^2$ and the definition of $\varphi_1$.
 
Consequently, starting from $0\le \varphi_2< \varphi_{1}<1$ on $(0, +\infty)$, we show by induction that $0 \le \varphi_k <\varphi_{k-1}<1$ on $(0, +\infty)$ for every $k\ge 2$. One checks again by induction that $\varphi_k$ is continuous since by change of variable
\[
\varphi_{k}(t) =  \varphi_1(t)- \rho\int_0^t f^2_{\lambda}(s)\varsigma^2(t-s)\big( 1- \varphi_{k-1}(t-s)\big) \frac{ds}{\lambda^2 c}
\]
(since $\varsigma^2$ is bounded and continuous owing to Assumption~$E_{\lambda}$).

By the first Dini Lemma, it follows that $\varphi_k\downarrow  \varphi_{\infty}\!\in {\cal C}(\R_+,\R)$ uniformly on compact intervals  with $\varphi_{\infty}(0)=1$. In particular $\varphi_{\infty}$ satisfies the functional equation 
\[
\varphi_{\infty}(t)  =  R_{\lambda}^2(t)+ \rho\int_0^t f^2_{\lambda}(t-s)\varsigma^2(s) \varphi_{\infty}(s) \frac{ds}{\lambda^2 c}.
\]

Let $\displaystyle \ell_{\infty} := \limsup_{t\to +\infty} \varphi_{\infty}(t) \!\in [0,1]$. For every $\varepsilon >0$ there exists $A_{\varepsilon} >0 $ such that for $t\ge A_{\varepsilon} $, $\varphi_{\infty}(t) \le \ell_{\infty}+\varepsilon$. Then 
\begin{align*}
\int_0^t f^2_{\lambda}(t-s)\varsigma^2(s) \varphi_{\infty}(s) \frac{ds}{\lambda^2 c} & \le \int_{A_{\varepsilon} }^t f^2_{\lambda}(t-s) \varsigma^2(s) (\ell+\varepsilon) ds+ \frac{\|\varsigma\|_{\infty}}{c\lambda^2}\int_{t- A_{\varepsilon}}^tf^2_{\lambda}(u)du.
\end{align*}
Consequently, as $f_{\lambda}\!\in L^2({\rm Leb}_1)$ and $\lim_{t\to+\infty}R^2_{\lambda}(t)=0$,
\[
\limsup_{t\to +\infty} \varphi_{\infty}(t) \le \rho(\ell_{\infty}+\varepsilon)
\]
so that $\ell_{\infty}\le \rho\, \ell_{\infty}$ which in turn implies $\ell_{\infty}=0$ since $\rho \!\in [0,1)$. \hfill$\Box$

\bigskip
\noindent{\bf Remark.} If $[\sigma]_{\rm Lip}^2 <\frac{\lambda ^2} {\|\varsigma_{\lambda,c}^2\|_{\infty}\int_0^{+\infty}f^2_{\lambda}(u)du}<1$ and $R_{\lambda}\!\in {\cal L}^2({\rm Leb}_1)$,  then one derives from Fubini-Tonelli's theorem that  
$$
\displaystyle \int_0^{+\infty} \varphi_{\infty}^2(s)ds\le \frac{\lambda^2}{\lambda^2-[\sigma]_{\rm Lip}^2\|\varsigma^2\|_{\infty}\int_0^{+\infty}f^2_{\lambda}(u)du}\int_0^{+\infty}R^2_{\lambda}(t)dt<+\infty.
$$

\subsection{Long run functional weak behaviour}
\begin{Theorem} \label{thm:longrun} Let $\lambda$, $c >0$. Assume the kernel $K$ and its $\lambda$-resolvent $R_{\lambda}$  satisfy
\begin{equation}\label{eq:flambda2beta}
  \int_0^{+\infty}f_{\lambda}(t)^{2\beta}dt <+\infty\quad \mbox{ for some }\beta>1
\end{equation}
and assume that  there exists $\vartheta\!\in (0,1]$ such that 
\begin{equation}\label{eq:flambda2theta}
  \int_0^{+\infty}(f_{\lambda}(t+\delta)-f_{\lambda}(t))^2 dt\le C\delta^{2\vartheta}.
\end{equation}
  Let $\lambda,c>0$, let $\mu_0\!\in \R$, and let $\sigma:\R\to \R$ be a Lipschitz continuous function. Assume $(E_{\lambda})$ is in force.
  Assume $f_{\lambda}\!\in L^2({\rm Leb}_1)$. Let    $(X_t)_{t\ge 0}$ be  the solution to the Volterra equation~\eqref{eq:Volterrameanrevert3} starting  from $X_0\!\in L^2(\P)$. 

\smallskip 
\noindent $(a)$ {\em $C$-tightness of the time shifted processes}. Assume furthermore that 
 \begin{equation}\label{eq:p&c}
X_0\in L^p(\P) \quad \mbox{ and }\quad  \left\{\begin{array}{lll}
 p=2 &\mbox{and } \;\;c<\frac{1}{[\sigma]^2_{\rm Lip}} & \mbox{ if }\;\; \vartheta\wedge  \frac{\beta-1}{2\beta}> \tfrac 12\\
  p> (\vartheta\wedge  \frac{\beta-1}{2\beta})^{-1} &\mbox{and } \; \;c<\frac{1}{4p[\sigma]^2_{\rm Lip}} & \mbox{ if }\; \;(\vartheta\wedge  \frac{\beta-1}{2\beta}) \le \tfrac 12.\end{array}\right.
 \end{equation}

Then,  the family of  shifted processes $(X_{t+u})_{u\ge 0}$  is  $C$-tight and uniformly  integrable, square uniformly integrable if $p>2$, as $t\to+\infty$. Thus,    for 
 any limiting distribution $P$  on the canonical space $\Omega_0:=C(\R_+, \R)$  equipped with the Borel $\sigma$-field induced by the sup-norm topology, the canonical process $Y_t(\omega)= \omega(t)$, $\omega \!\in \Omega_0$, has, for any small enough $\eta>0$, a   $\big((\vartheta\wedge \frac{\beta-1}{2\beta})-\frac 1p-\eta\big)$-H\"older pathwise continuous $P$-modification. 
 Moreover the shifted processes of two solutions $(X_t)_{t\ge 0}$  and $(X'_t)_{t\ge 0}$ are $L^2$-confluent in the sense that there exists a non-increasing function $\widetilde \varphi_\infty: \R_+\to [0,1]$   and $\lim_{t\to +\infty} \widetilde \varphi_{\infty} = 0$ such that
for every $0\le t_1<t_2< \cdots < t_{_N}<+\infty$
\[
{\cal W}_2\big([(X_{t+t_1}, \cdots, X_{t+t_{_N}})], [(X'_{t+t_1}, \cdots, X'_{t+t_{_N}})])\le \tilde \varphi_{\infty}(t)\to 0 \mbox{ as }t\to +\infty
\]
As  a consequence, the  functional weak limiting distributions  of $[X_{t+\cdot}]$ and $[X'_{t+\cdot}]$ are the same in the sense that,  if $[X_{t_n+\cdot}]\stackrel{(C)}{\longrightarrow} P$ for some subsequence $t_n \to +\infty$,
then  $[X'_{t_n+\cdot}]\stackrel{(C)_w}{\longrightarrow} P$ and  conversely.

%
%
%
%
\smallskip
\noindent $(b)$ {\em Functional weak long run behaviour}. Assume furthermore that the above solution $(X_t)_{t\ge 0}$ has a fake stationary regime of type~I starting from a random variable$X_0$ with mean $\frac{\mu_0}{\lambda}$ and variance $v_0$.
Then, for every $t_1,\, t_2 \ge 0$,  $t_1\le t_2$,
\begin{equation}\label{eq:Cinfty}
{\rm Cov}(X_{t+t_1}, X_{t+t_2})\stackrel{t\to+\infty}{\longrightarrow} C_{f_{\lambda}}(t_1,t_2):= \frac{v_0}{\int_0^{+\infty}f^2_{\lambda}(s)ds}\int_0^{+\infty}  f_{\lambda}(t_2-t_1+u)f_{\lambda}(u)du.
\end{equation}
Hence,    
 under any limiting distribution $P$, the canonical process   $Y$ is a  (weak) $L^2$-stationary processes with mean $\frac{\mu_0}{\lambda}$ and covariance function $C_{f_{\lambda}}(s,t)$, $s,\, t\ge 0$.

\smallskip 
\noindent $(c)$ When $\sigma(x)= \sigma>0$ is constant and $X_0$ has a normal distribution ${\cal N}(\frac{\mu_0}{\lambda}, v_0)$ , then  $(X_t)_{t\ge 0}$ satisfies
\[
X_{t+\cdot}\stackrel{(C)}{\longrightarrow } \Xi^{(f_{\lambda})} \quad \mbox{ as}\quad t\to +\infty,
\]
where   $\Xi^{(f_{\lambda})}$ is  the stationary Gaussian process with  mean $\frac{\mu_0}{\lambda}$ and covariance function $C_{f_{\lambda}}$ ($\stackrel{(C)}{\rightarrow }$ stands for functional weak convergence on $C(\R_+,\R)$ equipped with the topology of uniform convergence  on compact sets).
\end{Theorem} 

To be compared to the more precise result from~\cite{friesen2022volterra} for the Volterra  CIR Volterra process.

\bigskip
\noindent{ \bf Proof.}   $(a)$   It follows from~\eqref{eq:p&c} that either $p=2$  and $c<\frac{1}{[\sigma]^2_{\rm Lip}}$ or $p>2$ and $c<\frac{1}{[4p\sigma]^2_{\rm Lip}}$. Hence Proposition~\ref{prop:Momentctrl}  implies that $\sup_{t\ge 0}\Big \|X_t-\frac{\mu_0}{\lambda}\Big\|_p<+\infty$.   As a consequence, $\sigma$ having an at most affine growth as a Lipschitz function, we derive that $\sup_{t\ge 0} \|\sigma(X_t)\|_p<+\infty$ Now we are in position to establish  $C$-tightness by Kolmogorov criterion.
Let $s,\,t\ge 0$, $s\le t$. Let $p$ be given by~\eqref{eq:p&c}. 
\[
X_t-X_s = \big(R_{\lambda}(t)-R_{\lambda}(s)\big) \Big(X_0-\frac{\mu_0}{\lambda}\Big) + \bigg(\int_0^tf_{\lambda}(t-u)\varsigma(u)\sigma(X_u)dW_u- \int_0^sf_{\lambda}(s-u)\varsigma(u)\sigma(X_u)dW_u\bigg).
\]
 Let us denote by $(A)$ and $(B)$ the two terms of the sum on the right hand side of the above equality.
\begin{align*}
 \| A \|_p &= \Big\| X_0-\frac{\mu_0}{\lambda} \Big\|_p\int_s^t f_{\lambda}(u)du\\
& \le \Big\| X_0-\frac{\mu_0}{\lambda} \Big\|_p\left( \int_0^{+\infty}f_{\lambda}^{2\beta}(u)du\right)^{\frac{1}{2\beta}}|t-s|^{1-\frac{1}{2\beta}}\\
&=C_{X_0, \beta, f_{\lambda}} |t-s|^{1-\frac{1}{2\beta}}
\end{align*}
where $C_{X_0, \beta, f_{\lambda}}<+\infty$ owing  to~\eqref{eq:flambda2beta}.       On the other hand
\begin{align*}
 \| B \|_p &\le \left\|\int_s^t f_{\lambda}(t-u)\varsigma(u)\sigma(X_u)dW_u \right\|_p + \left\|\int_0^s \big(f_{\lambda}(t-u)-f_{\lambda}(s-u)\big)\varsigma(u)\sigma(X_u)dW_u \right\|_p. 
\end{align*}
Let us denote by $(B.1)$ and $(B.2)$ the two terms on  the right hand side. Combining the $L^p$-BDG and the generalized Minkowski inequalities yields for $(B.1)$ 
\begin{align*}
\|(B.1)\|_p &\le C^{BDG}_p \Big(\int_s^tf^2_{\lambda}(t-u) \varsigma^2(u)\|\sigma(X_u)\|_p^2 du\Big)^{1/2} \\
&\le C^{BDG}_p\sup_{u\ge 0} \|\sigma(X_u)\|_p \|\varsigma\|_{\infty}\Big(\int_s^tf^2_{\lambda}(t-u) du\Big)^{1/2} \\
&\le C^{BDG}_p\sup_{u\ge 0} \|\sigma(X_u)\|_p \|\varsigma\|_{\infty}\Big(\int_0^{+\infty}f^{2\beta}_{\lambda}(u) du\Big)^{1/2\beta} |t-s|^{\frac 12 (1-\frac{1}{\beta})}\\
&\le C_{p, \sigma, f_{\lambda}}|t-s|^{ \frac{\beta -1}{2\beta}}.
\end{align*}

As for $(B.1)$, one derives from~\eqref{eq:flambda2theta} that
\begin{align*}
\|(B.2)\|_p &\le C_{p, \sigma, f_{\lambda}}\left(\int_0^{+\infty} \big(f_{\lambda}(t-s+v)-f_{\lambda}(v)\big)^2du\right)^{1/2}\le C_{p, \sigma, f_{\lambda}}|t-s|^{\vartheta}.
\end{align*}

Finally, after noticing that $ \frac{\beta-1}{2\beta}\le 1-\frac{1}{2\beta} $ what precedes proves the existence of a real constant $C_{p,\lambda}>0$ such that 
\[
\E\, |X_t-X_s|^p \le C_{p,\lambda}|t-s|^{p(\vartheta\wedge \frac{\beta-1}{2\beta})}.
\]
Under our assumption~\eqref{eq:p&c}, one has $p(\vartheta\wedge \frac{\beta-1}{2\beta})>1$.  Then it follows from  Kolmogorov's $C$-tightness criterion (see~\cite[Theorem~2.1]{RevuzYor}), that the family of  shifted processes $X_{t+\cdot}$, $t\ge 0$, is $C$-tight with limiting distributions $P$ under which the canonical process has the announced H\ "older pathwise regularity. The confluence result follows from Proposition~\ref{prop:confluence} with $\widetilde \varphi_{\infty} (t)= \sup_{u\ge t} \varphi_{\infty}(u)$. 


\smallskip
 \noindent{\sc Step 5}. Let $X$ and $X'$ be two solutions of  Equation~\eqref{eq:Volterrameanrevert} starting from $X_0$ and $X'_0$ respectively, both square integrable. Using Proposition~\ref{prop:confluence}, we derive that for every $0\le t_1<t_2< \cdots < t_{_N}<+\infty$
\[
{\cal W}_2\big([(X_{t+t_1}, \cdots, X_{t+t_{_N}})], [(X'_{t+t_1}, \cdots, X'_{t+t_{_N}})])\to 0 \mbox{ as }t\to +\infty
\]
As  a consequence, the  weak limiting distributions  of $[X_{t+\cdot}]$ and $[X'_{t+\cdot}]$ are the same in the sense that,  if $[X_{t_n+\cdot}]\stackrel{(C)}{\longrightarrow} P$ for some subsequence $t_n \to +\infty$ (where $P$ is a probability measure on $C(\R_+, \R)$ equipped with the Borel $\sigma$-field induced by the sup-norm topology), then  $[X'_{t_n+\cdot}]\stackrel{(C)_w}{\longrightarrow} P$ and  conversely.

\smallskip
\noindent $(b)$ {\em Asymptotic weak $L^2$-stationarity}.
%
By an integration by part and using that $\varsigma(0)=0$ and $\varsigma$ is bounded, one derives that 
\[
L_{(\varsigma^2)'}(t)= t\, L_{\varsigma^2}(t) = \frac{2c \lambda^2L_{R_{\lambda}f_{\lambda}}}{L_{f^2_{\lambda}}}(t)
\]
where the second equality follows from~\eqref{eq:Laplacesigma}. As $f^2_{\lambda} \!\in {\cal L}^2({\rm Leb}_1)$ and $2L_{R_{\lambda}f_{\lambda}}(0) = - \int_0^{+\infty} (R^2_{\lambda})'(u)du = 1-0 = 1$, one has
\[
L_{(\varsigma^2)'}(t)\stackrel{0}{\sim} \frac{c\lambda^2}{\int_0^{+\infty}f^2_{\lambda}(s)ds}.
\]
By the Hardy-Littlewood Tauberian Theorem (see~\cite{FellerII}), we get that
\[
\varsigma^2(t) = \int_0^t (\varsigma^2)'(s)ds \stackrel{+\infty}{\sim}  \frac{c\lambda^2}{\int_0^{+\infty}f^2_{\lambda}(s)ds}.
\]

Now let us consider the asymptotic  covariance between $X_{t+t_1}$ and $X_{t+t_2}$, $0<t_1<t_2$ when $X_t$ starts for $X_0$ with mean $\frac{\mu_0}{\lambda}$,  $v_0$ and $\bar \sigma^2 = \E\, \sigma(X_t)^2$, $t\ge 0$.
\begin{align*}
{\rm Cov}(X_{t+t_1}, X_{t+t_2})& = {\rm Var}(X_0)R^2_{\lambda}(t)+ \frac{1}{\lambda^2}\E\int_0^{t+t_1} f_{\lambda}(t+t_2-s)f_{\lambda}(t+t_1-s)\varsigma^2(s)\E \,\sigma^2(X_s)ds\\
&=  {\rm Var}(X_0)R^2_{\lambda}(t)+ \frac{\bar \sigma^2}{\lambda^2}\E\int_0^{t+t_1} f_{\lambda}(t_2-t_1+u)f_{\lambda}(u)\varsigma^2(t+t_1-u)du.
\end{align*}
Note that $f_{\lambda}\!\in {\cal L}^2({\rm Leb}_1)$ by ~\eqref{eq:flambda2beta}  since $f_\lambda$ is a probability density. Hence $ f_{\lambda}(t_2-t_1+\cdot)f_{\lambda}\!\in {\cal L}^2({\rm Leb}_1)$ since $f_{\lambda}\!\in {\cal L}^2({\rm Leb}_1)$ and $\mbox{\bf 1}_{\{0\le u \le t+t_1\}}\varsigma^2(t+t_1-u)\to  \frac{c\lambda^2}{\int_0^{+\infty}f^2_{\lambda}(s)ds}$ for every $u\!\in \R_+$ as $t\to +\infty$,
\[
{\rm Cov}(X_{t+t_1}, X_{t+t_2})\stackrel{t\to+\infty}{\longrightarrow}  \frac{v_0}{\int_0^{+\infty}f^2_{\lambda}(s)ds}\int_0^{+\infty}  f_{\lambda}(t_2-t_1+u)f_{\lambda}(u)du 
\]
where we also used that $R^2_{\lambda}(t)\to 0$ as $t\to +\infty$.  Combining this with Sep~1 completes the proof of this claim.

\smallskip
\noindent $(c)$ This  follows first  from the fact that $(X_t)_{t\ge 0}$ is then a Gaussian process  so that its functional weak limiting distributions are Gaussian and secondly that a Gaussian process is entirely determined by its mean and covariance functions. \hfill$\Box$

\bigskip

\bigskip
\noindent {\bf Remarks.} $\bullet$ We would need $\|\sup_{t\ge 0}|X_t|\|_p<+\infty$ for some  $p>2$ to derive in $(a)$ relative compactness for (functional) ${\cal W}_2$-compactness (quadratic Wasserstein distance). 
 
\smallskip
\noindent $\bullet$ Be aware that at this stage we do not have uniqueness of the limit distributions since they are not characterized by their lean and covariance functions, except in Gaussian setting ($X$ is Gaussian process).

\medskip
\noindent {\bf Examples} $(a)$ When $\sigma(x)=\sigma>0$ and $X_0$ has a Gaussian distribution, the process $X$ is  Gaussian  and this proves (at least for finite dimensional weak convergence i.e. weak convergence of all marginals at any order).
\[
(X_{t+\cdot})\stackrel{{\cal D}_f}{\longrightarrow} X^{(\infty)} \quad \mbox{ as }t\to +\infty
\]
where $X^{(\infty)}$ is a Gaussian process with mean $\frac {\mu_0}{\lambda}$ and covariance function given by 
\[
{\rm Cov}(X^{(\infty)}_{t_1}, X^{(\infty)}_{t_2})= v_0 \int_0^{+\infty}  f_{\lambda}(t_2-t_1+u)f_{\lambda}(u)\frac{du}{\int_0^{+\infty}f^2_{\lambda}(v)dv} .
\]
\noindent $(b)$ Let  $\sigma$ be   given the Lipschitz continuous function with $[\sigma]_{\rm Lip} =\sqrt{\kappa_2}$ defined by~\eqref{eq:sigmafakeI&II} i.e. $\sigma(x) =  \sqrt{\kappa_0+\kappa_1(x-\frac{\mu_0}{\lambda})+\kappa_2(x-\frac{\mu_0}{\lambda})^2}$ with $\kappa_1^2-4\kappa_0\kappa_2\le 0$ (like in Proposition~\ref{prop:quasiZumbach}).  Let $c\!\in (0, 1/\kappa_2)$. If  $X_0\!\in L^2(\P)$ with $\E\, X_0 = \frac{\mu_0}{\lambda}$ and ${\rm Var}(X_0)= v_0= \frac{c\kappa_0}{1-c\kappa_2}>0$ then the  solution $(X_t)_{t\ge 0}$ to~\eqref{eq:Volterrameanrevert} has a fake stationary regime of type~I. If $\sigma$ is constant it becomes of type~II.

As for the long run behaviour described in claim~$(a)$ of the above theorem, it holds whenever   $X_0\!\in L^p(\P)$ such that $p$ and $c$ satisfy~\eqref{eq:p&c} (which depend on parameters of the kernel $K$) with $c$ given by Proposition~\ref{prop:quasiZumbach}. In particular the condition on $c$ reads en terms of variance $v_0<\frac{\kappa_0}{\kappa_2}\frac{4p}{4p-1}$ so that, if $X_0\!\in  \cap_{p>2}L^p(\P)$, the condition reads  $v_0 < \frac{\kappa_0}{\kappa_2}$.
%
%
%

\section{Back to $\alpha$-fractional kernels, $0<\alpha<1$}\label{sec:Back2ML}
This case is of special interest since it corresponds to the recent introduction of  fractional Volterra SDEs to devise ``rough models'' of stochastic volatility dynamics in Finance (see~\cite{JaissonR2016,BayerFG2016,GatheralJR2018,ElEuchR2018}) or long memory volatility models when $H\!\in (1/2,1)$ (see~\cite{ComteR1996, ComteR1998}). In these models Volterra equations with fractional kernels $K_\alpha(t)=\frac{ t^{\alpha-1}}{\Gamma(\alpha)}$, $\alpha= H+\frac 12$ with  $H\!\in (0,\frac 12)$ appear as more  tractable dynamics than solutions of SDEs involving  stochastic integrals with respect to true $H$-fractional Brownian motions ($H$ being the Hurst coefficient) which would require to call upon  ``high order" rough path theory.

Our aim in this section is to prove that for such  kernels $K_{\alpha}$, the results obtained in  former Sections \textcolor{red}{\dots}  apply. First that  the resolvent $R_{\alpha,\lambda}$  satisfies   the monotonicity assumption~\eqref{eq:hypoRlambda} for every  $\lambda>0$ and that   $f_{\alpha,\lambda}= -R_{\alpha,\lambda} $ exists and  is square integrable w.r.t. the Lebesgue   measure on $\R_+$, at least when  $\frac 12<\alpha<1$. The main result of this section is of course that he stabilizers $\varsigma_{\alpha, \lambda,c}$ exist (and share additional specific properties) so that the results established in Sections~\ref{sec:TowardFake} and~\ref{sec:FakeI-II}  in the case of interest $\sigma(t,x)= \sigma(t)$ (Gaussian setting) and $\sigma(t,x)= \varsigma(t)\sigma(x)$ apply (in particular their existence when $\sigma(x)$ is given by~\eqref{eq:sigmafakeI&II}).

To this end, it will be enough to study   the case $\lambda =c =1$ in most  situations. Let $R_{\alpha}=R_{\alpha,1}$  given by its expansion~--~denoted $e_\alpha$ in the literature  and known as {\em Mittag-Leffler function}~--~given by 
\[
R_\alpha(t)= e_\alpha(t) = \sum_{k\ge 0} (-1)^k \frac {t^{\alpha k}}{\Gamma(\alpha k +1)}, \quad t\ge 0.
\]
since then
\begin{equation}\label{eq:scaling}
R_{\alpha, \lambda}= R_\alpha(\lambda^{1/\alpha}\cdot) = e_\alpha(\lambda^{1/\alpha}\cdot)
\quad\mbox{ 
and }
\quad
f_{\alpha, \lambda}= -R'_{\alpha, \lambda}= \lambda^{1/\alpha}f_\alpha(\lambda^{1/\alpha}\cdot)  
\end{equation}
where
\[
f_\alpha:= f_{\alpha,1} = - e'_\alpha.
\]
The functions $R_\alpha$  (and subsequently $f_\alpha$) is {\em completely monotone} (CM) on the $(0, +\infty)$  in the sense that    $(-1)^n R_\alpha^{(n)} \ge 0$ at every order $n\ge 0$. 

This follows  from the fact that $e_{\alpha}$ is the Laplace transform  of a  non-negative Lebesgue  integrable function, namely
$$
e_{\alpha}(t)=L_{H_{\alpha}}(t) = \int_0^{+\infty}e^{-tu}H_{\alpha}(u)du
$$ 
 where $H_{\alpha} : (0, +\infty)\to \R_+$. This representation  was first established in~\cite{Pollard1948}. More recently, a  synthetic formula was found for $H_{\alpha}$ in~\cite{GorMain2000} (see (F.22) p.31, see also~\cite{Mainardi2014}). Let us provide more details. One starts from the Laplace transform of $e_{\alpha}$ which writes on $\mathbb{C}$,  
\[
L_{e_{\alpha}}(z)= \frac{z^{\alpha-1}}{z^{\alpha}+1},\quad z\in \mathbb{C}, \;\Re(z)>1
\]
where $z^{\alpha}:= |z|^{\alpha}e^{i\alpha{\rm arg(z)}}$, $-\pi < {\rm arg}(z) <\pi$.  The above identity is based on  the inverse Laplace transform (Bromwich-Mellin formula) 
\[
e_\alpha(t)= \frac{1}{2i\pi} \int_{s=1-i\cdot\infty}^{s=1+i\cdot\infty}e^{st}\frac{s^{\alpha-1}}{s^{\alpha}+1}ds
\]
which finally yields a closed   form for $H_\alpha$ given by 
\begin{equation}\label{eq:Halpha}
\forall\,  u\!\in (0,+\infty)\; H_{\alpha}(u)=- \frac{1}{2\pi} \cdot 2\,\Im {\rm m}\Big( \frac{z^{\alpha-1}}{z^{\alpha}+1}\Big)_{|z= ue^{i\pi}}  = \frac{\sin(\alpha\pi)}{\pi}\frac{u^{\alpha-1}}{u^{2\alpha}+2u^{\alpha}\cos(\pi \alpha)+1}>0
\end{equation}
since $0<\alpha <1$ implies $\sin(\alpha\pi)>0$ and 
\begin{equation}\label{eq:ineqHalpha}
u^{2\alpha}+2u^{\alpha}\cos(\pi \alpha)+1 \ge 1-\cos^2(\alpha \pi)= \sin^2(\alpha \pi)>0
\end{equation}

Moreover, as $H_\alpha$ is continuous on $(0,+\infty)$, 
$$
H_\alpha(u)\stackrel{0}{\sim}  u^{\alpha-1}\frac{\sin(\pi \alpha)}{\pi}
\quad \mbox{ and }\quad H_\alpha(u)\stackrel{+\infty}{\sim} \frac{\sin(\pi \alpha)}{\pi u^{\alpha+1}},
$$
 it is clear that $H_\alpha \!\in {\cal L}_{\R_+}^1({\rm Leb}_1)$ and that both functions
\[
u\mapsto uH_\alpha(u) \quad \mbox{ and }\quad u\mapsto u^{\alpha+1}H_\alpha(u) \mbox{ are bounded on $\R_+$.}
\]
Thus, for every $\varepsilon >0$, $\int_0^{+\infty} e^{-\varepsilon u}uH_\alpha (u)du <+\infty$ so that $e_\alpha$ is differentiable  on $(0, +\infty)$ with 
\begin{equation}\label{eq:e'alpha}
e'_\alpha(t) = -\int_0^{+\infty} e^{-tu} uH_\alpha(u) du< 0, \quad t>0.
\end{equation}
\begin{Proposition}Let $\lambda>0$ and let $\alpha\!\in (0,1)$.

\smallskip
\noindent $(a)$  The $\lambda$-resolvent $R_{\alpha, \lambda}$  satisfies $R_{\alpha, \lambda}(0)=1$ and 
is completely monotonic, hence infinitely differentiable on $(0, +\infty)$. Moreover, $R_{\alpha,\lambda}$ decreases  to $0$, $R_{\alpha,\lambda}\!\in {\cal L}^{r}({\rm Leb}_1)$ for every $r>\frac{1}{\alpha}$ and $f_{\alpha, \lambda}= -R'_{\alpha, \lambda}$ is a completely monotonic function  (hence convex), decreasing to $0$  and satisfying 
\[
\forall\, t>0, \quad f_{\alpha, \lambda} (t):= -R'_{\alpha, \lambda} (t) = \lambda^{\frac{1}{\alpha}}\int_0^{+\infty} e^{-\lambda^{\frac{1}{\alpha}}tu}uH_{\alpha}(u)du>0.
\]

\noindent $(b)$ Moreover, if $\alpha\!\in (\frac 12,1)$, $f_{\alpha,\lambda} \!\in {\cal L}^{2\beta}({\rm Leb}_1)$ for every $\beta \!\in \big(0, \frac{1}{2(1-\alpha)}\big)$  for every $\vartheta\!\in\big(0,\alpha-\frac 12\big)$,  there exists a real constant $C_{\vartheta, \lambda}>0$ such that 
\[
\forall\, \delta>0, \quad \left[\int_0^{+\infty} \big(f_{\alpha,\lambda}(t+\delta)-f_{\alpha,\lambda}(t) \big)^2 \right]^{1/2}\le C_{\vartheta,\lambda}\delta^{\vartheta}.
\]
  \end{Proposition}

\bigskip
\noindent {\bf Proof.} $(a)$ follows when $\lambda =1$ from the fact  that  $R_{\alpha}= e_{\alpha}$,   hence decreasing  to $0$, completely monotonic      with $f_{\alpha}= -e'_\alpha >0$. The representation of $f_\alpha$ follows from~\eqref{eq:e'alpha}  It follows  from~\eqref{eq:Halpha} and~\eqref{eq:ineqHalpha}    that 
\begin{equation}\label{eq:Halphabound}
H_\alpha(u) \le \frac{u^{\alpha-1}\sin(\pi \alpha)}{\pi \sin^2(\pi \alpha)} =  \frac{u^{\alpha-1}}{\pi \sin(\pi \alpha)}.
\end{equation}
Hence, for every $t\ge 0$,  
\[
R_{\alpha}(t) = e_{\alpha}(t)\le \frac{1}{\pi \sin(\pi \alpha)}\int_{0}^{+\infty}e^{-tu}u^{\alpha-1}du= \frac{\Gamma(\alpha)}{\pi \sin(\pi \alpha)}t^{-\alpha}
\]
so that $R_{\alpha}\!\in {\cal L}^{\gamma}({\rm Leb}_1)$ for every $\gamma>\frac{1}{\alpha}$. This extends to $R_{\lambda, \alpha}$ by scaling using~\eqref{eq:scaling}.
 
\smallskip
\noindent $(b)$ Let us prove the $L^{2\beta}$-integrability of $f_{\alpha, \lambda}$. Once again we may assume that $\lambda =1$ by scaling. noted that $f_{\alpha,\lambda}= \lambda^{1/\alpha} f_{\alpha}(\lambda^{1/\alpha} \cdot)$ so that $\int_0^{+\infty} f_{\alpha, \lambda}^{2\beta}(t)dt= \lambda^{\frac{2\beta-1}{\alpha}}\int_0^{+\infty} f_{\alpha}^{2\beta}(t)dt$, it is clear that it is enough to prove that $ f_{\alpha}$ is ${\cal L}^{2\beta}$-integrable.

\smallskip It follows from~\eqref{eq:Halphabound}  that for every $t>0$
 \[
f_{\alpha}(t) = - e'_\alpha(t) \le \frac{1}{\pi \sin(\pi \alpha)} \int_0^{+\infty} e^{-tu}u^{\alpha}  du= \frac{\Gamma(\alpha+1)}{t^{\alpha+1}\pi \sin(\pi \alpha)}.
\]

On the other hand 
\[
f_{\alpha}(t)\stackrel{0}{\sim} 
\frac{t^{\alpha-1}}{\Gamma(\alpha)}  
\]
owing~(\footnote{Another argument is to note that $uH_{\alpha}(u)  \stackrel{+\infty}{\sim} \frac{\sin(\pi \alpha)}{\pi}u^{-\alpha}$ so  that it follows from Karamata's theorem (see~\cite{BiGoTe1989}) that 
\[
f_{\alpha}(t)\stackrel{0}{\sim} \Gamma(1-\alpha) \frac{\sin(\pi \alpha)}{\pi}.
\]
The two formulas coincide since $\Gamma(\alpha)\Gamma(1-\alpha) =\frac{\pi}{\sin(\pi \alpha)}$.}) 
to~\eqref{eq:flambda}. Combining these two remarks implies the existence of  a real constant $C_\alpha>0 $ such that for $t\!\in (0,1]$, 
\[
 f_{\alpha}(t)\le C_{\alpha}\Big(\frac {1}{t^{1-\alpha}}\wedge \frac{1}{t^{\alpha+1}}\Big).
\]

As $t\mapsto \Big(\frac {1}{t^{1-\alpha}}\wedge \frac{1}{t^{\alpha+1}}\Big)\!\in {\cal L}^{2\beta}({\rm Leb}_1)$ for any $\beta\!\in \big( 1, \frac{1}{2(1-\alpha)}\big)$  since $\alpha \!\in (\frac 12, 1)$, $f_{\alpha}\!\in {\cal L}^{2\beta}({\rm Leb}_1)$. 

\smallskip
\noindent Another consequence is that, for every $t\ge 1$, 
$$
R_{\alpha}(t) =e_{\alpha}(t) = \int_t^{+\infty}f_{\alpha}(s)ds\le C'_{\alpha}t^{-\alpha}
$$ 
so that $R_{\alpha}\!\in {\cal L}^2({\rm Leb}_1)$.

\smallskip
As for the ${\cal L}^2({\rm Leb}_1)$-$\vartheta$-H\"older continuity of $f_{\alpha, \lambda}$, one may again assume w.l.g. that $\lambda =1$. Let $\delta >0$. One has
\[
f_{\alpha }(t+\delta)-f_{\alpha}(t) = \int_0^{+\infty} e^{-tu}(1-e^{-\delta u} ) uH_{\alpha}(u)du.
\]
As $0\le 1-e^{-v} \le (1-e^{-v})^{\vartheta}$, for every $v\ge 0$, since $\vartheta\!\in (0,1)$
\begin{align*}
\int_0^{+\infty} \big(f_{\alpha }(t+\delta)-f_{\alpha}(t) \big)^2 dt & \le \int_{(0, +\infty)^2}\!\!\! (uv)^{1+\vartheta}H_{\alpha}(u)H_{\alpha}(v)\int_0^t \!\!e^{-t(u+v)} dt\,du\,dv\, \delta^{2\vartheta} \\
&=   \int_{(0, +\infty)^2}\!\!\! \frac{(uv)^{1+\vartheta}}{u+v}H_{\alpha}(u)H_{\alpha}(v )\,du\,dv\, \delta^{2\vartheta} \\
& \le  \tfrac 12 \int_{(0, +\infty)^2}\!\!\!  (uv)^{\frac 12 +\vartheta} H_{\alpha}(u)H_{\alpha}(v )\,du\,dv\, \delta^{2\vartheta} \\
& = \tfrac 12 \left[\int_{(0, +\infty)}\!\!\!  u^{\frac 12 +\vartheta} H_{\alpha}(u)\,du\right]^2\, \delta^{2\vartheta},
\end{align*}
where we used  Fubini-Tonelli's theorem in the first   line to interwind the order of integration and the elementary inequality $\sqrt{uv} \le \frac 12 (u+v)$ when $u,\, v\ge 0$ in the penultimate line.  Now, we derive form~\eqref{eq:Halpha} that
$$
 H_\alpha(u)\stackrel{0}{\sim} \frac{\sin(\pi \alpha)}{\pi} u^{\alpha-1}
\quad \mbox{ and }\quad H_\alpha(u)\stackrel{+\infty}{\sim} \frac{\sin(\pi \alpha)}{\pi}u^{-(\alpha+1)},
$$
Consequently
$$
u^{\frac 12+\vartheta} H_\alpha(u)\stackrel{0}{\sim} \frac{\sin(\pi \alpha)}{\pi} u^{\alpha-\frac 12+\vartheta}
\quad \mbox{ and }\quad u^{\frac 12+\vartheta} H_\alpha(u)\stackrel{+\infty}{\sim} \frac{\sin(\pi \alpha)}{\pi} u^{-(\frac 12 +\alpha-\vartheta)},
$$
which implies that 
\[
\int_{(0, +\infty)}\!\!\!  u^{\frac 12 +\vartheta} H_{\alpha}(u)\,du<+\infty \quad \mbox{ iff }\quad \vartheta < \alpha-\tfrac 12.
\]
One concludes when $\lambda >0$ by scaling.
\hfill$\Box$

\smallskip
This leads to the following theorem e.g. for volatility coefficients of the form~\eqref{eq:sigmafakeI&II}.   First note that if if $\alpha \!\in \big(\frac 12, 1\big)$, then $\frac{1}{2\alpha-1}>2$;

 \begin{Theorem}$(a)$ Let  $\alpha \!\in \big(\frac 12, 1\big)$, let $K(t) = K_{\alpha}(t)= t^{\alpha-1}$, $t>0$, let $\sigma$ be given by~\eqref{eq:sigmafakeI&II},  let $p>\frac{1}{2\alpha-1}>2$,  let $c\!\in \big( 0,\frac{2\alpha-1}{4\kappa_2}\big)$ and let   $X_0\!\in L^p(\P)$  for some $p>\frac{1}{2\alpha-1}$ such that  $\E\, X_0= \frac{\mu_0}{\lambda}$ and ${\rm Var}(X_0)=\frac{c\kappa_0}{1-c\kappa _2}$. Then,

\smallskip
$(i)$  the solution $(X_t)_{t\ge 0}$ solution to the Volterra equation~\eqref{eq:Volterrameanrevert} starting from $X_0$ has a fake stationary regime of type~I,   

\smallskip
$(ii)$   the family of shifted processes $X_{t+\cdot}$, $t\ge 0$, is $C$-tight as $t\to +\infty$ and its (functional)  limiting distributions are all $L^2$-stationary processes with covariance function $C_{\infty}$ given by~\eqref{eq:Cinfty}.

\smallskip
$(iii)$  For every $X'_0\!\in L^2(\P)$ the solution to the equation~\eqref{eq:Volterrameanrevert} starting from $X'_0$ satisfies $\|X'_t-X_t\|_2\to 0$ as $t\to +\infty$.

\end{Theorem} 

\medskip \noindent {\bf Proof.} We may assume w.l.g. that $ \frac{1}{2\alpha-1}<p<\frac{1}{4c\kappa_2}$.
If $0<\vartheta < \alpha -\frac 12$ and $1< \beta<  \frac{1}{2(1-\alpha)}$, then $\displaystyle \lim_{\vartheta \to \alpha- \frac 12, \beta\to \frac{1}{2(1-\alpha)}}p(\big(\tfrac{\beta-1}{2\beta} )\wedge \vartheta\big )= \alpha\wedge(2\alpha-1)= p (2\alpha-1)>1$. Then Theorem~\ref{thm:longrun} applies.\hfill$\Box$

\subsection{Existence and computation of the function $\varsigma^2_{\a,\lambda,c}$ solution to Equation~\eqref{eq:Laplacesigma} when $\alpha\!\in (\frac 12,1)$}\label{sec:sigma2rough}

The aim of this  section is two-fold: we want simultaneously to prove the existence  of the  functions $\varsigma^2_{\alpha,\lambda,c}$ as the product of a monomial $t^{1-\alpha}$ with   power series in $t^{k\alpha}$ with an infinite convergence radius  and provide a way to compute it since since the coefficients  of the afore mentioned  expansion will be be easily computable by induction. At this stage the notation  $\varsigma^2_{\alpha, \lambda,c}$ should be understood in a synthetic  way, not as a square. Proving that  $\varsigma^2_{\alpha, \lambda,c}$ is non-negative is postponed to the next section and  relies on totally different arguments.

 To this end we rely on the Laplace version~\eqref{eq:Laplacesigma} of the equation satisfied by $\varsigma^2:=\varsigma^2_{\alpha,\lambda,c}$ (temporary notation), namely
 \begin{equation*}
\forall\, t>0, \quad  t\,L_{f^2_\lambda}(t).L_{\varsigma^2}(t)=  2\,c\lambda^2 L_{R_{\lambda}f_{\lambda}}(t).
\end{equation*}
 
Given the form of the kernel $K_{\alpha}(u) = \frac{u^{\alpha-1}}{\Gamma(\alpha)}$ and the expansion of the resolvents $R_{\lambda}$ (we drop the dependence in $\alpha$ for simplicity) and it derivative $-f_{\lambda}$, we check that
\[
R_{\lambda}f_{\lambda}(t)\stackrel{0}{\sim} \frac {\lambda t^{\alpha-1}}{\Gamma(\alpha)} \quad \mbox{ and}\quad f^2_{\lambda}(t)\stackrel{0}{\sim} \frac {\lambda^2 t^{2(\alpha-1)}}{\Gamma(\alpha)^2}
\]
so that -- at least heuristically~(\footnote{We use here heuristically a dual version of  Hardy--Littlewood's Tauberian theorem for Laplace transform, namely  $\varsigma^2(t)\stackrel{0}{\sim}Ct^{\gamma}$, $\gamma>-1$,  iff $L_{\varsigma^2}(t)\stackrel{+\infty}{\sim}C t^{-(\gamma+1)} \Gamma(\gamma+1).$})~--
\[
L_{R_{\lambda}f_{\lambda}}(t)\stackrel{+\infty}{\sim} \lambda t^{-\alpha} \quad \mbox{ and}\quad L_{f^2_{\lambda}}(t)\stackrel{+\infty}{\sim} \frac {\lambda^2\Gamma(2\alpha-1) t^{-(2\alpha-1)}}{\Gamma(\alpha)^2}.
\]
This implies that 
\[
L_{\varsigma^2}(t) \stackrel{+\infty}{\sim} 2\lambda \,c \frac{\Gamma(\alpha)^2}{\Gamma(2\alpha-1)} t^{-(2-\alpha)} 
\]
owing to Equation~\eqref{eq:Laplacesigma}.
This  in turn suggests that 
$$
\varsigma^2(t)\stackrel{0}{\sim}\frac{2\lambda c\Gamma(\alpha)^2}{\Gamma(2\alpha-1)\Gamma(2-\alpha)}\,t^{1-\alpha} \quad\mbox{ (so that $ \varsigma(0)=0$ since $\alpha <1$).}
$$
This suggests to search $\varsigma_{\a,\lambda,c}^2(t)$ of the form
\begin{equation}\label{eq:expectedsigma2}
\varsigma_{\a,\lambda,c}^2(t):= 2\,\lambda \, c\,t^{1-\alpha}\sum_{k\ge 0} (-1)^k c_k\lambda^k t^{\alpha k},
\end{equation}
with $c_0 = \frac{\Gamma(\alpha)^2}{\Gamma(2\alpha-1)\Gamma(2-\alpha)}$. Note such an identity implies   $\varsigma^2_{\a,\lambda, c}(0)=0$.

\medskip It is important to note at this stage that,  {\em $\a$ being fixed, all functions   $\varsigma_{\a,\lambda,c}^2$ from~\eqref{eq:expectedsigma2} are generated   by the same function}
\begin{equation}\label{eq:varsigmaref}
\varsigma_{\a}^2(t):= 2\,t^{1-\alpha}\sum_{k\ge 0} (-1)^k c_k t^{\alpha k}
\end{equation}
(with in mind that the $c_k$ depend on $\a$) by the scaling formula
\begin{equation}\label{eq:varsigmascaling}
\varsigma^2_{\a,\lambda, c}(t) = c \lambda^{2-\frac{1}{\alpha}}\varsigma_\alpha^2(\lambda^{1/\alpha} t).
\end{equation}
So we may assume in what follows that $c=\lambda =1$. 

Let us establish a  recurrence formula satisfied by the coefficients $c_k$ which makes possible the computation of the function $\varsigma^2_{\a}$. To this end, it is convenient (but not mandatory) to switch to Laplace transforms. 

First we can rewrite  the expansions~\eqref{eq:e_alpha} and~\eqref{eq:flambda} that define $R:=R_{\alpha,1}$ and $f= f_{\alpha,1}$ (temporary notation) as
\[
R (t)  =\sum_{k\ge 0} (-1)^k a_k t^{\alpha k}\quad \mbox{ with } \quad a_k = \frac{1}{\Gamma(\alpha k+1)}, \; k\ge 0,
\]
and 
\[
f (t) = t^{\alpha-1} \sum_{k\ge 0} (-1)^k b_k  t^{\alpha k} \quad \mbox{ with }b_k = \frac{1}{\Gamma(\alpha(k+1))}, \; k\ge 0.
\]
Then, using that $L_{u\mapsto u^{\gamma}}(t)= t^{-(\gamma+1)}\Gamma(\gamma+1)$, 
\[
L_{R f }(t) = \lambda t^{-\alpha}\sum_{k\ge 0} (-1)^k (a*b)_k  t^{-\alpha k}\Gamma(\alpha(k+1))
\]
and
\[
L_{f^2 }(t) =    t^{-2\alpha+1}\sum_{k\ge 0} (-1)^k (b^{*2})_k  t^{-\alpha k}\Gamma(\alpha(k+2)-1)
\]
where, for two sequences of real numbers $(u_k)_{k\ge 0}$ and $(v_k)_{k\ge 0}$, $(u*v)_k = \sum_{\ell=0}^ku_\ell v_{k-\ell}$. Set
\[
\widetilde b_k = ( b^{*2})_k \Gamma(\alpha(k+2) -1) \quad \mbox{ and}\quad \widetilde c_k = c_k \Gamma(\alpha(k-1)+2), \; k\ge 0. 
\]
Then if $\varsigma^2(t)$ has the expected form~\eqref{eq:expectedsigma2}, one has
\[
L_{\varsigma^2}(t) =2\,  t^{\alpha-2} \sum_{k\ge 0}(-1)^k \widetilde c_k  t^{-\alpha k} 
\]
so that 
Equation~\eqref{eq:Laplacesigma} reads, by identification of the coefficients of the expansions on both sides of the equation,
\[
\forall\, k\ge 0, \quad (\widetilde b*\widetilde c)_k  =  (a*b)_k \Gamma(\alpha(k+1)).
\]
Elementary computations confirm  that 
\[
c_0=\frac{\Gamma(\alpha)^2}{\Gamma(2\alpha-1)\Gamma(2-\alpha)}
\]
and show that, for every $k\ge 1$,
\small
\begin{equation}\label{eq:c0}
c_k = \frac{\Gamma(\alpha)^2}{ \Gamma(\alpha(k-1)+2)\Gamma(2\alpha-1)}\left[ \Gamma(\alpha(k+1))(a*b)_k-\sum_{\ell=1}^k\Gamma(\alpha(\ell+2)-1) \Gamma(\alpha(k-\ell-1)+2) (b^{*2})_{\ell}  c_{k-\ell}  \right].
\end{equation}
\normalsize
Using the classical identities $\Gamma(a)\Gamma(b) = \Gamma(a+b) B(a,b)$, $a,\, b>0$ where $B(a,b)= \int_0^1 u^{a-1}(1-u)^{b-1}du$ and $\Gamma(a+1)= a\Gamma(a)$, one gets the alternative formulation for the $c_k$ which is much  more appropriate for numerical computations:  for every $k\ge 1$,

\small
\begin{equation}\label{eq:ck2}
\hskip-0,15cm c_k = \frac{\Gamma(\alpha)^2\Gamma(\alpha(k+1))}{\Gamma(2\alpha-1)\Gamma(\a k+2-\a)}\left[ (a*b)_k- \alpha(k+1)\sum_{\ell=1}^kB\big(\alpha(\ell+2)-1,\alpha(k-\ell-1)+2\big) (b^{*2})_{\ell}  c_{k-\ell}  \right].
\end{equation}
\normalsize 

Following classical techniques already used in~\cite{CalGraPag21} among many others, one shows that the convergence radius  $\rho_{\alpha}$  of  the  power series  $\sum_{k\ge 0}c_kt^{\alpha k}$ satisfies
$$
\rho_{\alpha} =\big( \liminf_n \big((c_n)^{1/n}\big)\big)^{-1/\alpha}=+\infty.
$$
 This is proved in Proposition~\ref{prop:Rzetainfty} in Appendix~A when $\alpha \!\in (\frac 12,1)$.
 
\subsection{Existence of $\varsigma_{\a,\lambda, c}$ (i.e. positiveness of   $\varsigma^2_{\a,\lambda, c}$) and more\dots} 
%
%

First note that we may assume that $\lambda=c=1$ by the scaling property~\eqref{eq:varsigmascaling}. We still denote by $R$ and $f$ the resolvent   $R_{\alpha,1}$ and the opposite of its derivative $f_{\alpha,1}$ respectively.  We know from~\eqref{eq:VolterraVarbter} that $\varsigma_{\alpha}^2$ satisfies the abstract  equation~\eqref{eq:abstraite} below. Hence, the  proposition below applied with with $\gamma= 1-\alpha$ and $\kappa= 2\lambda^{2-1/\alpha} c$ combined with the existence of $\varsigma^2_{\alpha, \lambda,c}$  given by~\eqref{eq:expectedsigma2} in the previous section   proves that, if $\alpha\!\in   (\frac 12, 1)$, \\
\centerline{\em  the stabilizer $\varsigma_{\alpha, \lambda, c}$ exists as an increasing concave function on $[0,+\infty)$,  null at $0$ and}
\[
\lim_{t\to+\infty} \varsigma_{\alpha, \lambda, c}=\frac{\sqrt{2\,c}\,\lambda}{\|f_{\alpha, \lambda}\|_{L^2({\rm Leb}_1)}}.
\]
\begin{Proposition}\label{prop:morevarsigma}  We consider  the equation
\begin{equation}\label{eq:abstraite}
\forall\, t\ge 0, \quad \kappa\big(1-R^2(t) \big) = f^2*g(t),
\end{equation} 
where $\kappa$ is a positive real constant, $R:\R_+\to \R_+$ and $f=-R'$ satisfy $R(0)=1$, $R>0$ decreases to $0$ and $f>0$ .

\noindent $(a)$ Then~\eqref{eq:abstraite} has at most one solution in $L^1_{loc}({\rm Leb}_1)$. 

\smallskip 
\noindent $(b)$ Assume furthermore  that $f$ is decreasing and~\eqref{eq:abstraite} has a continuous solution $g$ defined on $[0,+\infty)$   such that   there exists $\gamma\ge 0$, $\alpha>0$, and an analytic function $\tilde g :\C \to \C$ such that 
\[
\forall\, t\ge 0,\; g(t) = t^{\gamma} \tilde g(t^{\alpha})\quad \mbox{ and }\quad \tilde g(0)> 0.
\]
Then $g\ge 0$ on $\R_+$ so that the function $\sqrt{g}$ is well-defined on $\R_+$.

\smallskip 
\noindent $(c)$ Still assume $({\cal K})$ is in force, $f\!\in {\rm Leb}^2(\R_+)$ is non-increasing on $(0,+\infty)$. Assume $0<\gamma <1$ and $\alpha+\gamma\ge 1$ in the above decomposition. Then $g$ is concave, non-decreasing  and  converges  toward a finite $\frac{\kappa}{\int_0^{+\infty}f^2(u)du}$.
\end{Proposition}

\noindent {\bf Proof}. $(a)$  If $g_1$ and $g_2$ are two solutions then  then $f^2*\delta g=0$ in $L^1_{loc}(\R_+)$ with $\delta g =g_1-g_2$. As $L_{f^2}>0$ then $L_{\delta g}=0$ which implies $\delta g=0$ in $L^1_{loc}(\R_+)$.

\smallskip
\noindent  $(b)$ Assume there exists $\tilde t\!\in (\varepsilon_0,T)$  such that $g(\tilde t)<0$. Then
 $t_0 =\sup\{t<\tilde t: g(t)<0\}<+\infty$. As $\tilde g(0)>0$, $g>0$ at least on a small interval $(0, \varepsilon_0]$ for some $\varepsilon_0>0$ so that $t_0 \ge \varepsilon_0>0$.  By continuity of $g$ it is clear $g(t_0)=0$ and $g\ge 0$ on $[0,t_0]$.
 
 Assume now that $g\le 0$ on a small interval $[t_0, t_0+\eta_0]$ for some $\eta_0>0$. Then, for every $t\!\in (t_0, t_0+\eta_0]$,
\begin{align*}
f^2*g(t) -f^2*g(t_0) & = \int_{t_0}^t \underbrace{f^2(t-s)}_{>0}\underbrace{g(s)}_{\le 0}ds +  \int_{0}^{t_0}\big( \underbrace{f^2(t-s)- f^2(t_0-s)}_{\le 0}\big)\underbrace{g(s)}_{\ge 0}ds\le 0
\end{align*}
since, for the second integral,  $f^2$ is non-increasing. On the other hand $f^2*g(t) -f^2*g(t_0) = \kappa  (R^2(t_0) -R^2(t))\ge 0$ which yields a contradiction. Hence, for every large enough $n\ge 0$, there exists $t^+_n\!\in (t_0, t_0+\frac 1n]$ such that $g(t^+_n)>0$. On the other  hand by the very definition of $t_0$ there exists a sequence $t^-_n>t_0$, $n\ge 1$ such that $g(t^-_n)<0$. One then builds by induction a sequence $(\tau_n)_{n\ge 1}$ such that $g(\tau_{2n+1})<0$ and $g(\tau_{2n})>0$, $\tau_n \to t_0$ as $n\to +\infty$, $\tau_n >t_0$.  In turn this implies by the intermediate value theorem the existence of a sequence $(\tilde \tau_n)_{n\ge 1}$ such that $\tilde g(\tilde \tau^{\alpha}_n) = g(\tilde \tau_n )=0$, $\tau^{\alpha}_n >t^{\alpha}_0$ and $\tau^{\alpha}_n \to t^{\alpha}_0$ by the continuity of $g$.  As $\tilde g$ is analytic, it implies that $\tilde g$ is everywhere zero. Hence a contradiction since $\tilde g(0)>0$.

\smallskip
\noindent $(c)$ It is clear that $g$ is ${\cal C}^{\infty}$ on $(0, +\infty)$, continuous on $[0,+\infty)$  and  $g(0)=0$ since $\tilde g$ is   analytic. One then checks that
\[
-\kappa (R^2)'= f^2*g'.
\]
Note that  $ -\kappa (R^2)'= 2\kappa\, Rf$ is decreasing and negative since both $R$ and $f$ are positive and decreasing. In particular $-\kappa (R^2)''< 0$ on $(0,+\infty)$.
%
%

It follows from the decomposition of $g$ that for every t$>0$,
\begin{align*}
g''(t)&= \gamma(\gamma-1)t^{\gamma-2}\tilde g(t)+ 2\alpha \gamma t^{\alpha+\gamma-1}\tilde g'(t^{\alpha}) +\alpha(\alpha+\gamma-1) t^{\alpha +\gamma -2}\tilde g'(t^{\alpha})+\alpha^2t^{2\alpha+\gamma-2}\tilde g''(t^{\alpha})\\
& = \gamma(\gamma-1)t^{\gamma-2}\tilde g(t)+\big((2\alpha\gamma t+(\alpha+\gamma-1))\tilde g'(t^{\alpha})+\alpha^2t^{\alpha}\tilde g''(t^{\alpha})\big)t^{\alpha+\gamma-2}.
\end{align*}
As $\alpha>0$, $t^{\alpha+\gamma-2}= o(t^{\gamma-2})$ when $t\to 0$. Moreover $2\alpha\gamma t+(\alpha+\gamma-1))\tilde g'(t^{\alpha})+\alpha^2t^{\alpha}\tilde g''(t^{\alpha})\sim (\alpha+\gamma-1)\tilde g'(0)$ as $t\to 0$.
Hence $g''<0$ on an interval $(0, \varepsilon_1]$ for some $\varepsilon_1>0$ since $0<\gamma <1$. Now assume there exists $\tilde t_1>0$ such that  $g''(\tilde t_1)>0$ and set
\[
t_1 = \inf \{t>0: g''(t)>0\}.
\]
Let $t\!\in [0, t_1)$ and $\varepsilon>0$ small enough such that $t+\varepsilon\!\in [0, t_1)$. Then
\begin{align*}
-\kappa\frac{ ( R^2)'(t+\varepsilon)- (R^2)'(t)}{\varepsilon}- \int_0^tf^2(s)\Big( \frac{g'(t+\varepsilon-s)-g'(t-s)}{\varepsilon}\Big)ds&=\frac{1}{\varepsilon} \int_0^{\varepsilon}f^2(t+\varepsilon-s)g'(s)ds.
\end{align*}
One easily checks that the left hand side of the above equality converges as $\varepsilon\to 0$ toward 
\[
\underbrace {-\kappa  (R^2)''(t)}_{<0}-\underbrace{\int_0^t f^2(s)g''(s)ds}_{\le 0} <0
\]
As for the right hand side, we have  
\begin{align*}
\frac{1}{\varepsilon} \int_0^{\varepsilon}f^2(t+\varepsilon-s)g'(s)ds&\ge f^2(t+\varepsilon)\frac{1}{\varepsilon} \int_0^{\varepsilon}g'(s)ds\\
& =  f^2(t+\varepsilon)\frac{g(\varepsilon)}{\varepsilon}\\
& \to +\infty\quad \mbox{ as }\quad \varepsilon \to 0.
\end{align*}
Hence a contradiction.
Consequently $g''\le 0$ on $(0, +\infty)$. Hence $g$ is concave and non-negative. Which implies that $g$ is non-decreasing and has a limit $g_\infty\!\in(0,+\infty]$ as $t\to +\infty$.
%
%
If $g_\infty = +\infty$, then for every $A>0$   there exists $t_A$ such that for $t\ge t_A$, $g(t)\ge A$. Hence
\begin{align*}
f^2*g(t) &\ge \int_0^{t_A}  f^2(t-s)g(s)ds + A \int_{t_A}^t f^2(t-s)ds\\
& =  \int_0^{t_A}  f^2(t-s)g(s)ds + A \int_{0}^{t-t_A} f^2(u)du.
\end{align*}
Consequently, as $f^2*g(t)= \kappa (1-R^2(t)) \to  \kappa $ as $t\to +\infty$, 
\[
\kappa = 
 \lim_{t\to +\infty} f^2*g(t)\ge  A \int_{0}^{+\infty} f^2(u)du. 
\]
As $f\!\in L^2(\R_+)$, this yields a contradiction by letting $A$ go to infinity. Hence $g_\infty<+\infty$ and the same arguments show that 
$g_{\infty} \le \frac{\kappa}{\int_0^{+\infty}f^2(u)du}$ by considering $A = A_{\eta}=g_{\infty}-\eta$, $\eta>0$ by letting $\eta\to 0$ a posteriori. By the same reasoning one shows that $g_{\infty} \ge \frac{\kappa}{\int_0^{+\infty}f^2(u)du }$ which completes the proof.  \hfill$\Box$ 

\begin{figure}[h!]
\centering
\begin{tabular}{cc}
 \hskip-1cm \includegraphics[height = 5.5cm, width= 9cm]{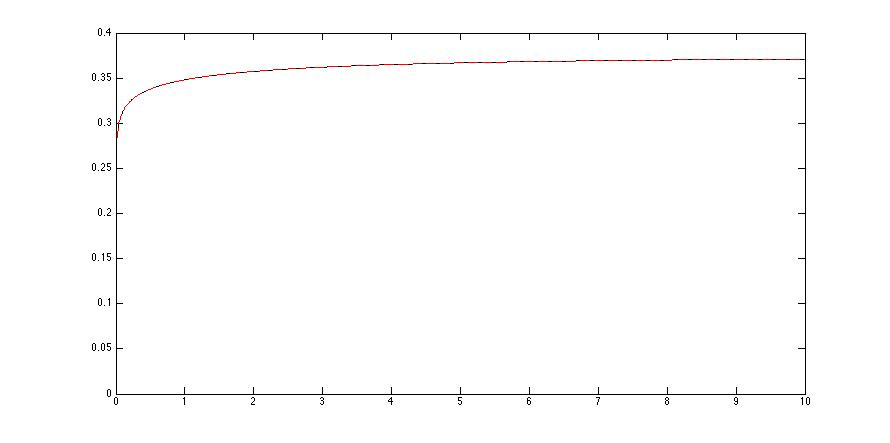}& \hskip-0,5cm\includegraphics[height = 5.5cm, width= 9cm]{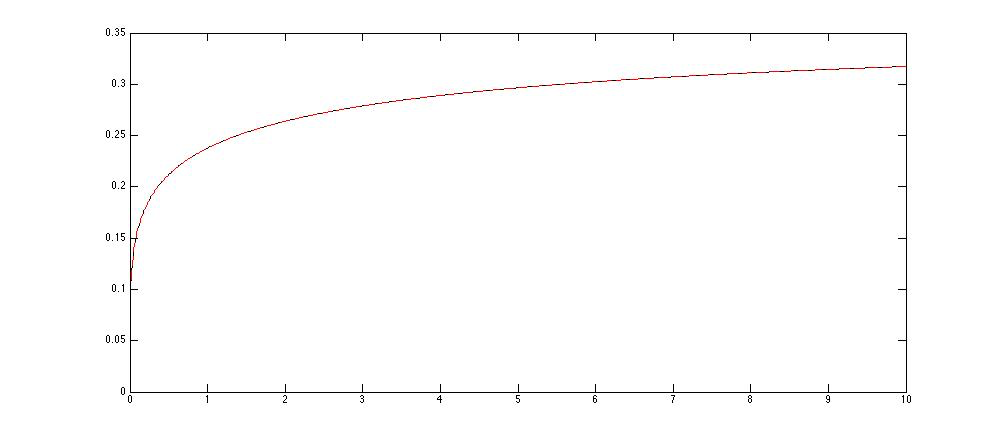}
\end{tabular}
\caption[]{\em Graph of $t\mapsto \varsigma_{\a,\lambda,c}(t)$ over time interval $[0,T]$, $T=10$. Left:  $H=0.4$, $\lambda = 0.2$, $c=0.3$. Right: $H=0.1$, $\lambda = 0.2$, $c=0.3$.}
\label{fig:1} 
 \end{figure}

\bigskip
\noindent $\blacktriangleright$ Typical examples of graphs is displayed in Figure~\ref{fig:1}  below 
 with parameters $H=0.4$ ($\alpha = 0.9$ (Left) and $H=0.1$ ($\alpha =0.6$ (Right), $\lambda = 0.2$, $v_0= 0.3$, $T=10$ (with $100$ steps per unit of time). These examples and many others strongly suggest that $\varsigma^2$  is non-negative in our framework (boundedness follows from the continuity of $\varsigma$ and its finite limit at infinity).
 
\subsection{A numerical illustration: (stabilized) quadratic rough Heston volatility }
We consider the scaled Volterra equation~\eqref{eq:Volterrameanrevert} with  kernel $K= K_{\alpha}$, mean-reverting coefficient $\lambda >0$,  long run mean $\frac{\mu_0}{\lambda}$ and  diffusion coefficient $\sigma(t,x)$ given by
\[
\sigma(t,x )= \varsigma(t) \sigma(x) \quad \mbox{with}\quad \sigma(x) = \sqrt{\kappa_0 +\kappa_2\big(x-\frac{\mu_0}{\lambda}\big)^2}, \quad \kappa_0, \, \kappa_2 >0, 
\]
where $X_0$ satisfies $\E\, X_0 = \frac{\mu_0}{\lambda}$, ${\rm Var}(X_0) = v_0>0$ and the stabilizer $\varsigma= \varsigma_{\alpha, \lambda,c}$ is solution to~\eqref{eq:VolterraVarbter} with $c =\frac{v_0}{\kappa_0+v_0\kappa_2}$. Then $\sigma(t,x)$ is $\sqrt{\kappa_2}$-Lipschitz in $x$ uniformly in $t\!\in [0,T]$ since $\varsigma $ is bounded and $c\kappa_2=  \frac{v_0\kappa_2}{\kappa_0+v_0\kappa_2}<1$ since $\kappa_0>0$. This is a special case of  Proposition~\ref{prop:quasiZumbach} where $\kappa_1 =0$ which is in accordance with the quadratic rough volatility model~in\cite{GaJuRo2020}.

\smallskip
We want to illustrate the fact that the variance of the  resulting Volterra process has constant variance. 

We introduce an Euler scheme with step $\frac Tn$ of the semi-integrated form~\eqref{eq:Volterrameanrevert2} of the  equation  (with $\mu(s)= \mu_0$ and noting that $\int_0^t f_{\lambda}(s)ds = 1-R_{\lambda}(t)$) which reads for every $k=0,\ldots,n$,  
\[
\bar X_{t_{k}} = \frac{\mu_0}{\lambda} + \big(X_0 -\frac{\mu_0}{\lambda} \big)R_{\lambda}(t_k)+\frac{1}{\lambda}\sum_{\ell=1}^{k}f_{\lambda}(t_k-t_{\ell-1})\varsigma(t_{\ell-1}) \sigma(\bar X_{t_{\ell-1}})(W_{t_{\ell}}-W_{t_{\ell-1}})
\]
where $t_k = \frac{kT}{n}$, $k=0,\ldots,n$. 

We set numerical values of the parameters of the Volterra equation to 
\[
H=0.4,\; \mu_0=2,\quad v_0= 0.09,\quad \lambda = 0.2, \quad \kappa_0 = 0.25, \quad \kappa_2 = 0.384
\]
so that  $c  \simeq 0.3163$ and  $c\kappa_2 \simeq 0.1215$ so that $4c \kappa_2 \simeq 0.4860< 0.8= 2(\frac 12 +0.4)-1$. Finally we specified $X_0$ as $X_0\sim {\cal N}\big(\tfrac{\mu_0}{\lambda};v_0\big)$.
 
\smallskip We set the number of time steps at $n= 1000$ and the horizon time at $T=1$. We performed a  Monte Carlo simulation of successive size $M= 
100\, 000$.

\begin{figure}[h!]
\centering
\includegraphics[height = 5.5cm, width= 15.5cm]{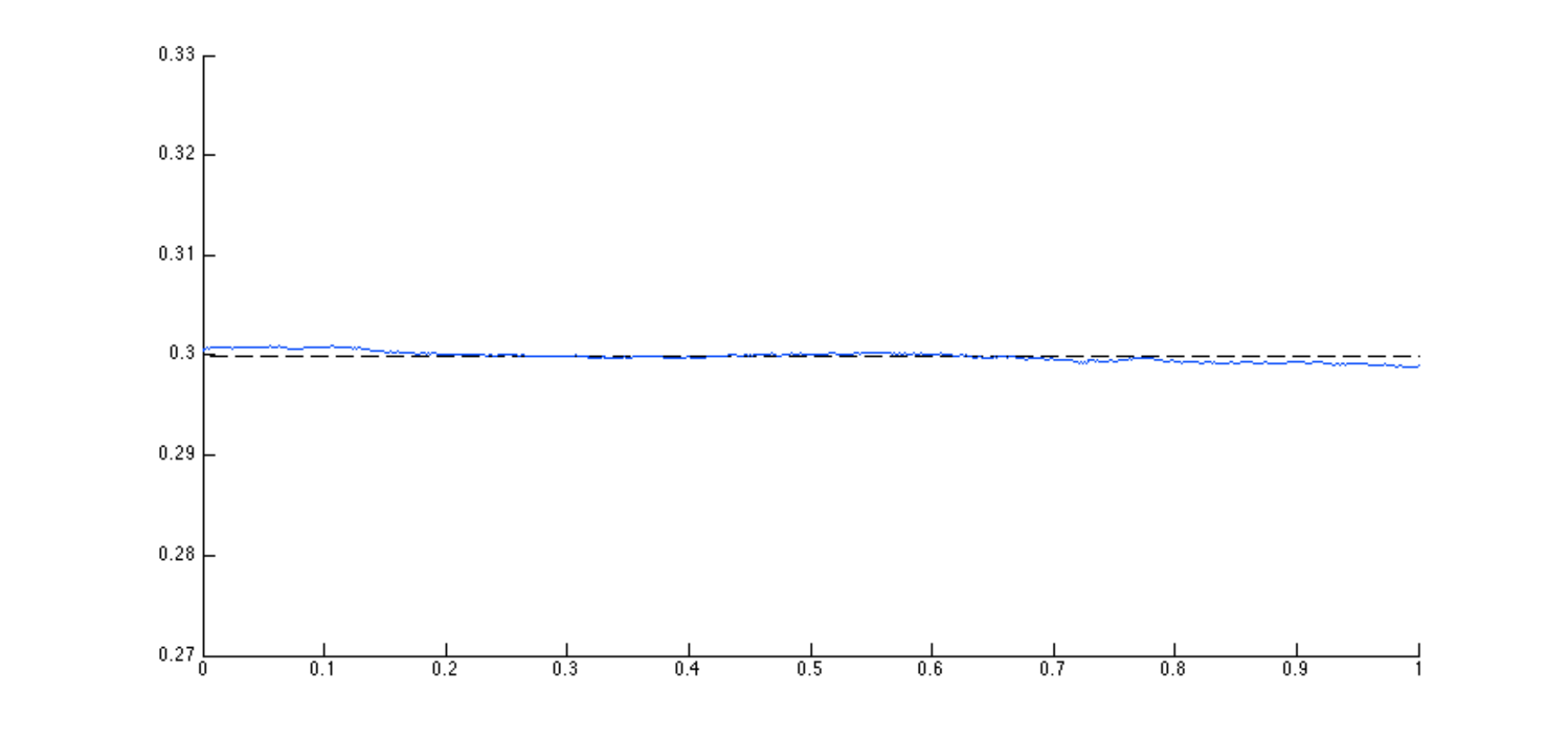}
 \caption{\em Graph of \textcolor{blue}{$t_k\longmapsto \sqrt{\overline{\rm Var}(t_k, M)}$} over the time interval $[0,T]$, $T=1$, $H=0.4$, $\mu_0=2$, $\lambda = 0.2$ , $v_0= 0.09$ and {\rm StdDev}$(X_0)= 0.3$ (- - -), $\kappa_0 = 0.25$, $\kappa_2 = 0.384$.
Number of steps: $n=1\,000$, Simulation size:  $M=100\,000$. }
\label{fig:fig3} 
 \end{figure}

\smallskip We depicted in Figure~\ref{fig:fig3} the empirical variances of the process at each discretization instant, namely 
$$
t_k \longmapsto \overline{\rm Var}(t_k, M)
$$
 of this simulation where  
 $$
 \overline{\rm Var}(t_k, M) = \frac{1}{M} \sum_{m=1}^M \big(\bar X^{(m)}_{t_k}-\frac{\mu_0}{\lambda}\Big)^2, \quad k=0,\ldots,n,
 $$
where $X^{(m)}$, $m=1,\ldots,M$  are i.i.d. copies of the Euler scheme. (We are aware that this  is not the  usual unbiased estimation of the variance in a Monte Carlo simulation).
%
 
 Note that, in beyond the clear stabilization effect,  one observes a significant numerical instability in spite of the fact that we considered $H=0.4$ (and not $H=0.1$) in this simulation.
%
%

\section{Provisional remarks} At this stage the  existence of non-trivial  fake stationary regimes when $\sigma$ is not given by~\eqref{eq:sigmafakeI&II} remains an open  question subject to ongoing research. By trivial fake stationary regimes we essentially mean true stationary regimes of homogeneous Brownian diffusion processes corresponding to the case of a constant kernels or for Volterra SDEs  (with non-constant kernels)  the cases where  $\sigma(\mu_0/\lambda)= 0$ so the mean  $\frac{\mu_0}{\lambda}$ is  common zeros of the affine drift and $\sigma$ (such is the case  of the fractional  CIR). however,  various numerical tests (carried out with kernel $K_{\alpha}$, $\frac 12 <\alpha<1$) suggest that  fake stationary regimes, at least  of type~I, exist for a wide class of volatility coefficients. 

\bigskip
\small

\noindent {\bf Acknowledgement:}  I thank N. Fournier, A. Pannier and M. Rosenbaum for helpful discussions.
\bibliographystyle{plain} 
\bibliography{propconvex}

\begin{thebibliography}{10}

\bibitem{alfonsi2023}
Aur\'elien Alfonsi.
\newblock Nonnegativity preserving convolution kernels. application to
  stochastic volterra equations in closed convex domains and their
  approximation, 2023.

\bibitem{BayerFG2016}
Christian Bayer, Peter Friz, and Jim Gatheral.
\newblock Pricing under rough volatility.
\newblock {\em Quant. Finance}, 16(6):887--904, 2016.

\bibitem{BiGoTe1989}
N.~H. Bingham, C.~M. Goldie, and J.~L. Teugels.
\newblock {\em Regular variation}, volume~27 of {\em Encyclopedia of
  Mathematics and its Applications}.
\newblock Cambridge University Press, Cambridge, 1989.

\bibitem{CalGraPag21}
Giorgia Callegaro, Martino Grasselli, and Gilles Pag\`es.
\newblock Fast hybrid schemes for fractional {R}iccati equations (rough is not
  so tough).
\newblock {\em Math. Oper. Res.}, 46(1):221--254, 2021.

\bibitem{CarlenK1991}
Eric Carlen and Paul Kr\'{e}e.
\newblock {$L^p$} estimates on iterated stochastic integrals.
\newblock {\em Ann. Probab.}, 19(1):354--368, 1991.

\bibitem{ComteR1996}
F.~Comte and E.~Renault.
\newblock Long memory continuous time models.
\newblock {\em J. Econometrics}, 73(1):101--149, 1996.

\bibitem{ComteR1998}
Fabienne Comte and Eric Renault.
\newblock Long memory in continuous-time stochastic volatility models.
\newblock {\em Math. Finance}, 8(4):291--323, 1998.

\bibitem{ElEuchR2018}
Omar El~Euch and Mathieu Rosenbaum.
\newblock Perfect hedging in rough {Heston} models.
\newblock {\em Ann. Appl. Probab.}, 28(6):3813--3856, 2018.

\bibitem{FellerII}
William Feller.
\newblock {\em An introduction to probability theory and its applications.
  {V}ol. {II}}.
\newblock John Wiley \& Sons, Inc., New York-London-Sydney, second edition,
  1971.

\bibitem{friesen2022volterra}
Martin Friesen and Peng Jin.
\newblock Volterra square-root process: Stationarity and regularity of the law,
  2022.

\bibitem{GatheralJR2018}
Jim Gatheral, Thibault Jaisson, and Mathieu Rosenbaum.
\newblock Volatility is rough.
\newblock {\em Quant. Finance}, 18(6):933--949, 2018.

\bibitem{GaJuRo2020}
Jim {Gatheral}, Paul {Jusselin}, and Mathieu {Rosenbaum}.
\newblock {The quadratic rough Heston model and the joint S\&P 500/VIX smile
  calibration problem}.
\newblock {\em arXiv e-prints}, page arXiv:2001.01789, January 2020.

\bibitem{GorMain1997}
Rudolf Gorenflo and Francesco Mainardi.
\newblock Fractional calculus: integral and differential equations of
  fractional order.
\newblock In {\em Fractals and fractional calculus in continuum mechanics
  ({U}dine, 1996)}, volume 378 of {\em CISM Courses and Lect.}, pages 223--276.
  Springer, Vienna, 1997.

\bibitem{Jaber2017}
Eduardo~Abi Jaber, Martin Larsson, and Sergio Pulido.
\newblock Affine volterra processes.
\newblock {\em The Annals of Applied Probability}, 2017.

\bibitem{Jacquieretal2022}
Antoine {Jacquier}, Alexandre {Pannier}, and Konstantinos {Spiliopoulos}.
\newblock {On the ergodic behaviour of affine Volterra processes}.
\newblock {\em arXiv e-prints}, page arXiv:2204.05270, April 2022.

\bibitem{JaissonR2016}
Thibault Jaisson and Mathieu Rosenbaum.
\newblock Rough fractional diffusions as scaling limits of nearly unstable
  heavy tailed {H}awkes processes.
\newblock {\em Ann. Appl. Probab.}, 26(5):2860--2882, 2016.

\bibitem{JouPag22}
Benjamin {Jourdain} and Gilles {Pag{\`e}s}.
\newblock {Convex ordering for stochastic Volterra equations and their Euler
  schemes}.
\newblock {\em arXiv e-prints}, pages arXiv:2211.10186, to appear in {\em Fin.
  \& Stoch.}, November 2022.

\bibitem{LemMonPag2022}
Vincent Lemaire, Thibaut Montes, and Gilles Pag\`es.
\newblock Stationary {H}eston model: calibration and pricing of exotics using
  product recursive quantization.
\newblock {\em Quant. Finance}, 22(4):611--629, 2022.

\bibitem{Mainardi2014}
Francesco Mainardi.
\newblock On some properties of the {M}ittag-{L}effler function
  {$E_\alpha(-t^\alpha)$}, completely monotone for {$t>0$} with {$0<\alpha<1$}.
\newblock {\em Discrete Contin. Dyn. Syst. Ser. B}, 19(7):2267--2278, 2014.

\bibitem{GorMain2000}
Francesco Mainardi and Rudolf Gorenflo.
\newblock Fractional calculus: special functions and applications.
\newblock In {\em Advanced special functions and applications ({M}elfi, 1999)},
  volume~1 of {\em Proc. Melfi Sch. Adv. Top. Math. Phys.}, pages 165--188.
  Aracne, Rome, 2000.

\bibitem{PagPan2009}
Gilles Pag\`es and Fabien Panloup.
\newblock Approximation of the distribution of a stationary {M}arkov process
  with application to option pricing.
\newblock {\em Bernoulli}, 15(1):146--177, 2009.

\bibitem{Pollard1948}
Harry Pollard.
\newblock The completely monotonic character of the {M}ittag-{L}effler function
  {$E_a(-x)$}.
\newblock {\em Bull. Amer. Math. Soc.}, 54:1115--1116, 1948.

\bibitem{RevuzYor}
Daniel Revuz and Marc Yor.
\newblock {\em Continuous martingales and {B}rownian motion}, volume 293 of
  {\em Grundlehren der mathematischen Wissenschaften [Fundamental Principles of
  Mathematical Sciences]}.
\newblock Springer-Verlag, Berlin, third edition, 1999.

\bibitem{RiTaYa2020}
Alexandre Richard, Xiaolu Tan, and Fan Yang.
\newblock Discrete-time simulation of stochastic {V}olterra equations.
\newblock {\em Stochastic Process. Appl.}, 141:109--138, 2021.

\bibitem{ZhangXi2010}
Xicheng Zhang.
\newblock Stochastic {V}olterra equations in {B}anach spaces and stochastic
  partial differential equation.
\newblock {\em J. Funct. Anal.}, 258(4):1361--1425, 2010.

\end{thebibliography}
\normalsize 

\appendix
\small
\section{Existence and computation  of  $\zeta^2_\a$ and $\zeta_\a$}

The aim of this additional section is to prove that the power series involved in the definition~\eqref{eq:varsigmaref} of the function $\zeta_\a^2(t) = 2 t^{1-\a}\sum_{k\ge 0} c_k t^{k}$  has an infinite convergence radius. Positiveness follows form~Proposition~\ref{prop:morevarsigma}. We entirely rely on the definition by induction~\eqref{eq:ck2} of the coefficients $c_k$ and the definitions for the sequences $a_k= \frac{1}{\Gamma(\a k+1)}$ and $b_k= \frac{1}{\Gamma(\a (k+1))}$, $k\ge 0$. By the triangle inequality we get the bound
\small
\begin{equation}\label{eq:ineqck}
\hskip-0,5cm |c_k| \le  \frac{\Gamma(\alpha)^2\Gamma(\alpha(k+1))}{\Gamma(\a k+2-\a)\Gamma(2\alpha-1)}\left[ (a*b)_k+ \alpha(k+1)\sum_{\ell=1}^kB\big(\alpha(\ell+2)-1,\alpha(k-\ell-1)+2\big) (b^{*2})_{\ell}  |c_{k-\ell}|  \right].
\end{equation}
\normalsize 

\medskip
We will extensively use the standard identities, for $a$, $b>0$: $\Gamma(a+1)=a\Gamma(a)$, $B(a,b):= \frac{\Gamma(a)\Gamma(b)}{\Gamma(a+b)}$.

\begin{Proposition}\label{prop:Rzetainfty} For every $\a \!\in (1/2,1)$, the convergence radius of the  power series $\sum_{k\ge 0} c_kt^k$ is infinite. To be more precise, there exists $K\ge \frac{1}{\Gamma(2-\a)}$ and $A\ge 2^{\a}$ such that 
\[
\forall\, k\ge 0,\qquad |c_k| \le K \frac{A^k}{\Gamma(\a k +2-\a)}.
\]
As a consequence the expansion~\eqref{eq:expectedsigma2} holds for every $t\!\in\R_+$ (in fact $\R$).
\end{Proposition}
\begin{Lemma}\label{lem:ubounds}Let $\a  \in (1/2,1)$. For every $k\ge 1$,

\smallskip
\begin{itemize}
\item $\displaystyle (a*b)_k \le \frac{2^{\a k} }{\Gamma(\a(k+1))}\big(1+(k+1) (1+\log k)\big)$.
\item $\displaystyle b^{*2}_k \le \frac{\a 2^{\a k}  }{\Gamma(\a(k+2))} (k+1)^2$.
\end{itemize}
\end{Lemma}

\noindent {\bf Proof}.  Using the identity linking the Beta function and the $\Gamma$ function, we have for every $k\ge 1$
\begin{align*}
(a*b)_k &= \sum_{\ell=0}^k \frac{1}{\Gamma(\a \ell)\Gamma(\a (k-\ell+1))}\\
&= \frac{1}{\Gamma(\a(k+1))}\bigg(1+  \sum_{\ell =1}^k \frac{1}{\ell} \frac{1}{B(\a \ell, \a (k-\ell+1))} \bigg).
\end{align*}

\smallskip
\noindent If $\ell=1$ or $k$, $B(\a \ell, \a (k-\ell+1)) = B (\a, \a k) = \int_0^1 u^{\a -1} (1-u)^{\a k -1} du \ge \frac{1}{\a k}$. 

\smallskip
\noindent If $2\le  \ell\le k-1$
\begin{align*}
B(\a \ell, \a (k-\ell+1)) &= \int_0^{1/2} u^{\a \ell-1} (1-u)^{\a(k-\ell+1)}du  = \int_0^{1/2} \ldots du + \int_{1/2}^1 \ldots du \\
&\ge \frac{1}{2^{\a (k+1)-1}} \Big( \frac{1}{\a \ell }+\frac{1}{\a (k-\ell+1)}\Big)\ge  \frac{1}{\a2^{\a (k+1)}} \frac{8}{k+1}\\
&\ge \frac{1}{\a2^{\a (k+1)}} \frac{4}{k+1} =  \frac{1}{\a2^{\a (k+1)-2}} \frac{1}{k+1}.
\end{align*}
The last downgraded inequality is also satisfied by when $\ell=1,k$.  Hence
\begin{align*}
(a*b)_k &\le  \frac{1}{\Gamma(\a(k+1))}\bigg(1 +\a 2^{\a(k+1)-2}(k+1)\sum_{\ell=1}^k \frac{1}{\ell}    \bigg) \\
 	    & \le  \frac{1}{\Gamma(\a(k+1))}\bigg(1 +\a 2^{\a(k+1)-2}(k+1)\sum_{\ell=1}^k\frac{1}{\ell}    \bigg)\\
	    & \le \frac{2^{\a k}}{\Gamma(\a(k+1))} \big(1+(k+1)(1+\log k)\big)
\end{align*}
since $\a 2^{\a-1}<1$ when $\a \!\in (0,1)$.

For $k=0$, $c_0 =  \frac{\Gamma(\a)^2\Gamma(2-\a)}{\Gamma(2\a-1)}\le \frac{K}{\Gamma(2-\a)}$.  Now let $k\ge 1$.

\smallskip
Let us deal now with $b^{*2}_k$ following the same principle.   
\[
b^{*2}_k = \frac{1}{\Gamma(\a(k+2))}\sum_{\ell=0}^k \frac{1}{B(\a (\ell+1), \a(k-\ell+1)}.
\]
If $\ell=0,k$,
\[
B(\a (\ell+1), \a(k-\ell+1)= B(\a, \a (k+1))=\int_0^1u^{\a -1} (1-u)^{\a(k+)-1}du\ge \frac{1}{\a(k+1)}
\]
and, if $1\le \ell\le k-1$, 
\[
B(\a (\ell+1), \a(k-\ell+1)\ge \frac{k+2}{\a 2^{\a(k+2)-1}(\ell+1)(k-\ell+1)} \ge  \frac{4}{\a 2^{\a(k+2)-3}(k+2)} .
\]
Hence, for every $k\ge 1$,
\begin{align*}
b^{*2}_k & \le \frac{1}{\Gamma(\a(k+2))}\bigg( 2\a (k+1+ \a \, 2^{\a -3}\sum_{\ell=1}^{k-1} 2^{\a k}(k+2)\bigg)\\
&= \frac{\a 2^{\a k}(k+1)^2}{\Gamma(\a(k+2))}\bigg(\underbrace{ \frac{2}{(k+1)2^{\a k}}+ 2^{\a -3} \frac{(k+1)(k+2)}{(k+1)^2}}_{=: \varphi_\a(k)} \bigg).
\end{align*}
One checks that $\sup_{k\ge1} \varphi_\a(k) <2^{-\frac 12}+ 2^{-2}<1$  for every $\a \!\in (1/2,1)$ which finally entails that, for every $k\ge 1$,
\[
\hskip 6cm b^{*2}_k \le \frac{\a 2^{\a k}}{\Gamma(\a(k+2))}(k+1)^2.\hskip 6cm\Box
\]

\bigskip
\noindent {\bf Proof of Proposition~\ref{prop:Rzetainfty}}. Let us prove now by induction that there exists $A>2^{\a}$, $K>|c_0| \Gamma(2-\a)$ such that 
\[
\forall\, k\ge 0, \quad |c_k| \le K \frac{A^k}{\Gamma(\a k+2-\a)}.
\]
Assume $c_\ell$  satisfies this inequality or every $\ell=0,\ldots, k-1$. Then, for every $\ell=1, \ldots, k$, 
\[
b^{*2}_\ell |c_{k-\ell}| \le K \frac{\a 2^{\a \ell}(\ell+1)^2 A^{k-\ell}}{\Gamma(\a(\ell+2))\Gamma(\a(k-\ell)+2-\a)}
\]
so that, noting that  $\Gamma(\a(k-\ell)+2-\a)= \Gamma(\a (k-\ell-1)+2)$,
\begin{align*} B(\a(\ell+2)-1,\a(k-\ell-1)+2))b^{*2}_\ell |c_{k-\ell}|  &\le   K \frac{\a 2^{\a \ell}(\ell+1)^2 A^{k-\ell}}{\a (\ell+2)\Gamma(\a(k+1)+1)}\\
&\le K \frac{2^{\a \ell}(\ell+1) A^{k-\ell}}{(k+1)\Gamma(\a(k+1)) } .
\end{align*}
Inserting this bound  into~\eqref{eq:ineqck} yields
\begin{align*}
|c_k|\le  \frac{\Gamma(\alpha)^2 }{\Gamma(2\alpha-1)\Gamma(\a k+2-\a)}\bigg(\frac{(a*b)_k}{\a(k+1)}+  K \,A^k\sum_{\ell=1}^k \Big(\frac{2^{\a}}{A} \Big)^{\ell}(\ell+1)  \bigg).
\end{align*}

Using the elementary inequality 
\[
\forall\, \rho \!\in (0, 1), \quad  \sum_{n\ge  1} n \rho^{n-1} \le \frac{1}{(1-\rho)^2}
\]
with $\rho= \rho(A): = 2^{\a}/A$, dividing the above inequality by $A^k$ and using the  upper bound for $(a*b)_k$  from Lemma~\ref{lem:ubounds} we get 
\[
\frac{|c_k|}{A^k}\le \frac{\Gamma(\alpha)^2 }{\Gamma(2\alpha-1)\Gamma(\a k+2-\a)}\bigg(\rho^k\a(1+\log k)+K \frac{1}{(1-\rho)^2}\bigg).
\]

Let $\varepsilon >0$ and let $A=A_{\a,\varepsilon}$ be large enough so that 
\[
\sup_{k\ge 1} \rho^k(k+1)(1+\log k)< \frac{\varepsilon}{\a}\quad\mbox{ and }\frac{1}{(1-\rho)^2}<1+\varepsilon
\]
then 
\begin{align*}
\frac{|c_k|}{A^k}& \le \frac{\Gamma(\alpha)^2 }{\Gamma(2\alpha-1)\Gamma(\a k+2-\a)}\bigg(\varepsilon+K(1+\varepsilon)\bigg) = \frac{K}{\Gamma(\a k+2-\a)} \frac{\Gamma(\alpha)^2 }{\Gamma(2\alpha-1)} \cdot \bigg(\frac{\varepsilon}{K}+  1+\varepsilon\bigg).
\end{align*}
To establish the heredity inequality we need to prove that 
\[
  \frac{\Gamma(\alpha)^2}{ \Gamma(2\alpha-1)} <1
\]
since then it is possible to choose $\varepsilon$ small enough and $K$ large enough so that $  \frac{\Gamma(\alpha)^2 }{\Gamma(2\alpha-1)} \cdot \big(\frac{\varepsilon}{K}+1+\varepsilon\big)
<1$.

Taking advantage of the $\log$-convexity of the $\Gamma$-function we have
\[
\log \Gamma(\a) \le \tfrac 12 \big( \log \Gamma(2\a-1) +\log \Gamma(1)\big) =   \tfrac12   \log \Gamma(2\a-1) 
\]
which implies the above inequality. 

Finally note that
\[
\varlimsup_k|c_k|^{\frac 1k}\le \varlimsup_k\bigg( \frac{A}{\Gamma(\a k +2-\a)}\bigg)^{1/k} =0
\]
since $\Gamma(\a k +2-\a)^{1/k}\sim e^{-\a} (\a k+2-\a)^{\a} \to +\infty \mbox{ as } k\to+\infty$ owing to Stirling's formula.\hfill$\Box$

\end{document}